\documentclass{article}
\usepackage{amsbsy,amssymb,amscd,amsfonts,latexsym,amstext,delarray,
amsmath}  \headsep=15pt
\def\pidr{ \pi_{1,DR}}
\def\piet{\pi_{1,\mbox{\'et}}}
\def\cP{ {\cal P} }
\def\Un{ \mbox{Un} }
\def\hpi{\hat{\pi}}

\def\picr{ \pi_{1,cr} }
\def\cX{{\cal X} }
\def\cD{ {\cal D}}

\def\End{\mbox{End}}
\def\bX{ \bar{X}}

\def\dim{{ \mbox{dim} }}
\def\Spec{{ \mbox{Spec} }}

\def\Hom{{ \mbox{Hom} }}

\def\Om{{ \Omega }}

\def\om{{ \omega  }}
\def\O{{ {\cal O} }}

\def\ra{{ \rightarrow }}
\def\da{{ \downarrow }}

\def\a{{ \alpha }}
\def\b{{ \beta }}
\def\g{{ \gamma }}

\def\Un{ \mbox{Un} }
\def\Vect{ \mbox{Vect} }
\def\hra{{ \hookrightarrow }}
\def\da{{ \downarrow }}

\def\C{{ \mathbb{C} }}

\def\bs{ \backslash}

\def\G{{ \Gamma }}
\def\Gal{{ \mbox{Gal} }}

\def\barF{ \bar{F}}

\def\bb{ \bar{b} }
\def\cE{ {\cal E}}

\def\Z{{ \mathbb{Z}}}

\def\Aut{ \mbox{Aut}}
\def\Isom{ \mbox{Isom} }

\def\uvb{ \mbox{unipotent vector bundle with connection}}
\newtheorem{thm}{Theorem}

\newtheorem{lem}{Lemma}
\newtheorem{prop}{Proposition}

\newtheorem{conj}{Conjecture}

\newtheorem{obs}{Observation}
\usepackage{amsmath}
\usepackage{amsfonts}
\usepackage{amscd}
\usepackage{amssymb}
\def\Q{\mathbb{Q}}

\def\invlim{\varprojlim}
\def\dirlim{\varinjlim}
\def\cV{ {\cal V} }
\def\cD{ {\cal D} }
\def\P{ {\bf P}}

\def\cT{ {\cal T}}
\def\jetloc{ j^{et}_{loc}}

\def\EM{ $\mbox{Eilenberg-Maclane}$ }

\title{The unipotent Albanese map and Selmer varieties for curves }
\author{ Minhyong Kim}

\begin{document}
\maketitle
\begin{abstract}
We study the unipotent Albanese map that associates
the torsor of paths for p-adic fundamental groups to a point
on a hyperbolic curve. It is shown that
the map is very transcendental in nature,
while  standard conjectures about the structure of mixed
motives provide control over the image of the map.
As a consequence, conjectures of `Birch and Swinnerton-dyer type'
are connected to finiteness theorems of Faltings-Siegel type.
\end{abstract}

\begin{em}
Dedicated to the memory of my teacher Serge Lang.
\end{em}

\medskip

In a letter to Faltings \cite{grothendieck} dated June, 1983, Grothendieck proposed
several striking conjectural connections between the arithmetic
geometry of `anabelian schemes' and their fundamental groups,
among which one finds issues of considerable interest to classical Diophantine
geometers. Here we will trouble the reader with
  a careful formulation of just one
of them. Let $F$ be a number field and $$f:X \ra \Spec(F)$$ a smooth, compact,
hyperbolic curve over $F$. After the choice of an algebraic closure
$$y:\Spec(\bar{F}) \ra \Spec(F)$$ and a  base point
$$x: \Spec (\bar{F}) \ra X$$ such that $f(x)=y$, we get an exact sequence
of fundamental groups:
$$
0\ra \hat{\pi}_1(\bar{X}, x) \ra \hat{\pi}_1(X, x)\stackrel{f_*}{ \ra} \G\ra 0,$$
where $\G=\Gal(\bar{F}/F)$ is the Galois group of $\bar{F}$
over $F$ and $\bar{X}=X\otimes \bar{F}$ is the base change of $X$
to $\bar{F}$. Now, suppose we are given a section
$s\in X(F)$ of the map $X\ra \Spec(F)$, i.e., an $F$-rational point
of $X$. This induces a  map of fundamental groups
$$s_*:\G \ra \hat{\pi}_1(X,x')$$ where $x'=s(y)$. Choosing an \'etale path
$p$ from $x'$ to $x$ determines an isomorphism
$$c_p:\hat{\pi}_1(X,x') \simeq \hat{\pi}_1(X,x), \ \ \ l \mapsto p\circ l\circ p^{-1},$$
 which is independent of
$p$ up to conjugacy. On the other hand, $p$ maps to
an element $\gamma \in \G$,  and it is straightforward to check
that $c_p\circ s_*\circ c_{\gamma^{-1}}$ is a continuous
splitting of
$f_*:\hat{\pi}_1(X,x)\ra \G$. This splitting is well-defined up
to the equivalence relation $\sim$ given by conjugacy of sections, where $\g \in \hat{\pi}_1(X,x)$ acts on
 a section $s$ as $s\mapsto c_{\g}\circ s \circ c_{f_*(\g)^{-1}}$.
We get thereby a map
$$s \mapsto [c_p\circ s_*\circ c_{\gamma^{-1}}]$$
from $X(F)$ to the set $Split (X))/\sim$ of
splittings
$Split (X)$ of the exact sequence  modulo the equivalence
relation $\sim$ given by conjugation. Grothendieck's {\em section conjecture}
states that this map is a bijection:
$$X(F) \simeq Split (X)/\sim $$
It seems that during the initial period of consideration, there
was an expectation that the section conjecture would
`directly imply' the Mordell conjecture. At present the
status of such an implication is unclear. Nevertheless,
this fact does not diminish the conceptual importance of the
section conjecture and its potential for broad
ramifications in
Diophantine geometry.  In the long run, one can hope that
establishing the correct link between Diophantine geometry
and homotopy theory will provide us with the framework
for a deeper understanding of Diophantine finiteness,
especially in relation to the analytic phenomenon of
hyperbolicity.

In this paper, we wish to make some preliminary comments
on the fundamental groups/Diophantine geometry
connection
from a somewhat different perspective which, needless to say,
does not approach the depth of the section conjecture.
In fact, in our investigation, the main tool is
 the {\em motivic fundamental group}, especially $p$-adic realizations,
rather than the pro-finite fundamental group. Nevertheless,
what emerges is a (murky) picture containing at least a few intriguing
points of mystery that are
 rather surprising in view of
the relative poverty of unipotent completions. That is to say, one is surprised by the
Diophantine depth of an invariant that is so close
to being linear.  It is the author's belief that
the {\em Selmer varieties} arising in this context, generalizing
$\Q_p$-Selmer groups of abelian varieties,
are objects of central interest.  Developing the formalism
for a systematic study of Selmer varieties is likely
to be crucial for continuing research along the
lines suggested in this paper. However, the non-abelian nature of the
construction presents a formidable collection of
obstacles that are at present beyond the author's power to surmount.
 In spite of this, it is hoped that
even a partial resolution of the problems can point us eventually towards the full
algebraic completion of the fundamental group, bringing a sort of
`motivic Simpson theory' to bear upon
the study of Diophantine sets.

We proceed then to a brief summary of the notions to be
discussed and
a statement of the results, omitting precise definitions
 for the purposes of this introduction. Although it will be clear that
at least part of the formalism is more general, we will focus
our attention on curves in this paper, right at the outset.
So we let $F$ be a number field and $X/F$ a smooth hyperbolic
curve, possibly non-compact. Let $R$ be the ring of $S-$integers
in $F$ for some finite set $S$ of primes. We assume that we
are given a smooth model
$$\begin{array}{ccc}
X & \hra & \cX \\
\da & & \da \\
\Spec(F) & \hra & \Spec(R)
\end{array}
$$
of $X$ and a compactification $\cX \hra \cX'$ relative to $R$
where $\cX'$  and the complement $\cD$ of $\cX$ in
$\cX'$ are also smooth over $R$. The Diophantine set of interest in
this situation is
$\cX(R)$, the $S$-integral points of $\cX$. The theorems of
Siegel and Faltings say that this set is finite.
After setting up the preliminary formalism, our goal will be to investigate these theorems
from a $\pi_1$-viewpoint.

As already mentioned, the $\pi_1$ we will be focusing on in this paper is a part of
the {\em motivic} $\pi_1$ defined by Deligne \cite{delignefg}.
In spite of the terminology, we do not rely on any general
theory of motives in our discussion of the fundamental group.
Rather, there will always be specific realizations related
by the standard collection of comparison maps.
Nevertheless, it is probably worth emphasizing that the
main idea underlying our approach is that of
a {\em motivic unipotent Albanese map}, of which the one defined by
Hain \cite{hain} is the Hodge realization.
That is, whenever a point $b\in \cX(R)$ is chosen as a basepoint,
 a unipotent motivic fundamental group
$U^{mot}:=\pi^{mot}_1(X,b)$
as well as motivic torsors $P^{mot}(y):= \pi^{mot}_1(X;b,y)$
of paths associated to other points $y$ should be
determined. The idea is that
if a suitable classifying space $D^{mot}$
for such torsors were to be constructed, the Albanese map would
merely associate to a point $y$, the class
$$[P^{mot}(y) ]\in D^{mot}.$$

Coming back to the concrete consideration of realizations, the analytic part of the machinery we need comes from the
{\em De Rham fundamental group} of $X_v:=X\otimes F_v$ that depends
only on the local arithmetic geometry. Here, $v$ is an archimedean
valuation of $F$ not contained in $S$. The definition of this
fundamental group requires the use of the
category $\Un (X_v)$ of unipotent vector bundles with
connections on $X_v$. Given any base point $b\in \cX(R_v)$,
we get a fiber functor $e_b: \Un (X) \ra \Vect_{F_v}$ to the
category of vector spaces over $F_v$, and
$$U^{DR}:=\pidr(X_v):= \Aut^{\otimes} (e_b)$$
in a rather obvious sense as functors on affine $F_v$-schemes.
The functor
$\pidr(X_v)$  ends up being representable by a pro-unipotent
pro-algebraic group. Natural quotients $$U^{DR}_n:=Z^{n}\backslash U^{DR}$$
via the descending central series of $U^{DR}$ correspond to restricting
$e_b$ to the subcategory generated by bundles having index of unipotency
$\leq n$. These quotients are  unipotent algebraic
groups over $F_v$. $U^{DR}$ is endowed with
a decreasing {\em Hodge filtration} by subgroups:
$$U^{DR} \supset \cdots \supset F^nU^{DR} \supset F^{n+1}U^{DR} \supset \cdots
\supset F^0U^{DR}$$
coming from a filtration of the coordinate ring by ideals.
There is a comparison isomorphism $$U^{DR} \simeq U^{cr}\otimes_KF_v
=\picr(Y_v, \bar{b})\otimes_KF_v$$
(where $K$ is the maximal absolutely unramified subfield of $F_v$)
with the crystalline fundamental group of the special fiber,
 defined using unipotent over-convergent iso-crystals. The utility of
this is that the Frobenius of the special
fiber comes to act naturally on the De Rham fundamental group.
Of crucial importance for us is the consideration of De Rham `path spaces':
$$P^{DR}(x):=\pidr(X_v;x,b):=\Isom^{\otimes} (e_b,e_x)$$
consisting of isomorphisms from the fiber functor $e_b$ to
$e_x$, $x\in X_v(R_v)$. These are pro-algebraic varieties
also endowed with Hodge filtrations
and crystalline Frobenii that are compatible with their structure
as right torsors for $U^{DR}$, and the point is to study how the
data vary with $x$. In fact, they turn out to be classified
by a natural `period space'
$$U^{DR}/F^0U^{DR}$$
giving us a higher De Rham unipotent Albanese map
$$j^{DR}: \cX(R_v) \ra U^{DR}/F^0U^{DR}$$
that maps $x$ to the class of $P^{DR}(x)$.
By passing to (finite-dimensional) quotients  $U^{DR}_n$, we also have
finite-level versions
$$j^{DR}_n: \cX(R_v) \ra U^{DR}_n/F^0U^{DR}_n$$
that fit into a compatible tower
$$\begin{array}{ccc}
|| & & \da   \\
\cX (R_v)& \ra & U^{DR}_{n+1}/F^0   \\
|| & & \da   \\
\cX (R_v)& \ra &   U^{DR}_n/F^0\\
|| &  & \da \\
 \vdots & \vdots & \vdots \\
|| &  & \da \\
\cX (R_v)& \ra & U^{DR}_2/F^0
\end{array}$$
At the very bottom, $U^{DR}_2=H_1^{DR}(X_v):=(H^1_{DR}(X_v))^*$
and $j^{DR}_2$ is nothing but the logarithm of the usual Albanese
map with respect to the basepoint $b$.

A comparison with the situation over $\C$  easily
yields the following:
\begin{thm} For each $n\geq 2$, the image of $j^{DR}_n$ is Zariski dense.
\end{thm}
This statement can be interpreted as linear-independence
for multiple polylogarithms of higher genus. It has been pointed out
 by a referee that this theorem
is also proved  by Faltings in a preprint \cite{faltings}. In
fact, both Faltings and Akio tamagawa had indicated to
the author earlier the possibility of proving the denseness
using trascendental methods in the genus zero case.

Unfortunately, this simple theorem is the only concrete result
to be reported on in this paper, the remaining parts being
an extended
commentary on what else might be expected when allowed considerable
optimism. Nevertheless, we proceed
to summarize here our observations.

The construction described thus far will have
to be compared with one involving pro-unipotent \'etale fundamental groups that  are associated to
the situation in both local and global settings.
That is, after a choice $\bar{F_v}$ of an algebraic
closure for $F_v$, we can associate to the
base-change $\bar{X}_v:= X_v\otimes \bar{F}_v$ the pro-unipotent
\'etale fundamental group $$U^{et}:=\piet (\bar{X}_v, b)$$ that classifies
unipotent lisse $\Q_p$-sheaves on $\bar{X}_v$. Furthermore,
to each $x$, we can also associate the space of unipotent
\'etale paths $P^{et}=\piet (\bar{X}_v; x,b)$ which is a torsor for
$U^{et}$. Both carry actions of
$G_v:=\Gal (\bar{F}_v/F_v)$ and the torsor structure is compatible with
the action. The torsors that are associated in this way
to integral points $x\in \cX(R_v)$  have the additional
property of being trivialized over a ring
$B_{cr}$ of $p$-adic periods. We classify these torsors
using a restricted Galois cohomology set
$H^1_f(G_v, U^{et})$. A classifying map defined exactly analogously to
the De Rham setting
provides a local \'etale Albanese map
$$\jetloc: \cX(R_v) \ra H^1_f(G_v, U^{et})$$
as well as finite-level versions
$$(\jetloc)_n: \cX(R_v) \ra H^1_f(G_v, U^{et}_n).$$
The connection to $j^{DR}$ goes through
a non-abelian extension of Fontaine's Dieudonn\'e functor,
 interpreted as an algebraic map
$$D: H^1_f(G_v, U^{et}_n) \ra U^{DR}_n/F^0$$
that fits into a commutative diagram
$$\begin{array}{ccc}
\cX(R_v) &\ra & U^{DR}_n/F^0 \\
 & \searrow & \uparrow \\
 & & H^1_f(G_v, U^{et}_n). \end{array} $$

The study of global points comes into the picture
when we take the base-point $b$ itself from $\cX(R)$,
the set of global integral points, and consider
the global $\Q_p$ pro-unipotent \'etale fundamental group
$\piet(\bar{X}, b)$ with the action of
$\G=\Gal(\bar{F}/F)$. A choice of an embedding
$\bar{F} \hra \bar{F_v}$ determines an inclusion
$G_v \hra \G$. Other than the consideration of the
Galois action, we have an isomorphism of pro-algebraic groups
$\piet(\bar{X}, b) \simeq \piet (\bar{X}_v, b)$,
so we will allow ourself a minimal abuse of notation
and denote the global fundamental group also by
$U^{et}$. Any other global point $y$ then determines
a torsor which is denoted $P^{et}(y)$ (again
introducing a bit of confusion with the local object).
 We denote by $T$ the set of primes consisting of
$S$ together with all primes dividing the residue characteristic
of $v$. The action of $G$ on $U^{et}$ factors through
$G_T$, the Galois group of the maximal subfield of $\barF$
unramified over the primes in $T$. Of crucial importance
for our purposes is the non-abelian cohomology set
$$H^1(G_T, U^{et})$$
and a natural subset
$$H^1_f(G_T, U^{et})\subset H^1(G_T, U^{et})$$
defined by `Selmer conditions' at the primes in $T$ (which
need not be too precisely determined in most situations).
All the cohomology sets thus far discussed can be interpreted
as the points of a pro-algebraic variety (over $\Q_p$), and
$H^1_f(G_T, U^{et})$ is the Selmer variety occurring in the title
of this paper. The finite-level versions
$H^1_f(G_T, U^{et}_n)$ and $H^1_f(G_v, U^{et}_n)$
are finite-dimensional algebraic varieties.
In this interpretation, the restriction
map and the Dieudonn\'e functor become algebraic maps of $\Q_p$-schemes,
the target of the latter being the Weil restriction
$Res^{F_v}_{\Q_p}(U^{DR}/F^0),$ which of course has the
property that
$$Res^{F_v}_{\Q_p}(U^{DR}/F^0)(\Q_p)=(U^{DR}/F^0)(F_v).$$
We have thus described the fundamental diagram
$$\begin{array}{ccccc}
\cX(R) &\ra &\cX(R_v)& \ra &U^{DR}_n/F^0 \\
\da & &\da & \nearrow & \\
H^1_f(G_T, U_n^{et})& \ra & H^1_f(G_v,U_n^{et}) &  &
\end{array}
$$
that provides for us the link between Diophantine geometry
and the theory of fundamental groups.
The reader familiar with the method of Chabauty \cite{chabauty}
will recognize here
a non-abelian lift of the diagram
$$\begin{array}{ccccc}
\cX(R) &\ra &\cX(R_v)& \ra &\mbox{Lie}(J )\otimes F_v\\
\da & &\da & \nearrow & \\
J(R)\otimes \Q_p& \ra &J(R_v)\otimes \Q_p &  &
\end{array}
$$
where $J$ is the Jacobian of $\cX$,
which is essentially the case $n=2$, that is, other than the replacement
of the Mordell-Weil group by the Selmer group.
The point is that the image of $\cX(R)$ inside
$U^{DR}_n/F^0$ thus ends up being contained inside the
image of $H^1_f(G_T, U_n^{et})$. The desired relation to Diophantine
finiteness is expressed by the
\begin{conj}
Let $F=\Q$ and $X$ be hyperbolic. Then
$$\dim (H^1_f(G_T, U_n^{et})) < \dim  (U^{DR}_n/F^0)$$
for some $n$ (possibly very large).
\end{conj}

Whenever one can find an $n$ for which this inequality of
dimensions is verified, finiteness of $\cX(R)$
follows, exactly as in Chabauty's argument.
The only difference is that the analytic functions that
occur in the proof are $p-$adic {\em iterated integrals} \cite{furusho}
rather than abelian integrals \cite{coleman}.

Unfortunately, this conjecture can be proved at present only
when the genus of $X$ is zero \cite{kim}
 and some special genus one situations (with a slight modification of the Selmer variety). However, there is
perhaps some interest in the tight connection we establish
between the desired statement and various conjectures
of `Birch and Swinnerton-dyer type.' For example,
our conjecture is implied by any one of:

(1) a certain fragment of the Bloch-Kato conjecture \cite{B-K};

(2) the Fontaine-Mazur conjecture \cite{F-M};

(3) Jannsen's conjecture on the vanishing of Galois cohomology
(when $X$ is affine) \cite{jannsen}.

\begin{flushleft}
In the case of CM elliptic curves (minus a point) it is
easily seen to follow also from the pseudo-nullity of a natural
Iwasawa module. All of these implications are rather easy once the formalism
is properly set up.
\end{flushleft}

It is amusing to note that  these conjectures
belong to what one might call `the structure theory of
mixed motives.' That is to say, the usual Diophantine
connection for them occurs through the theory of L-functions.
It is thus rather surprising
 that a non-linear (and non-trivial) phenomenon like Faltings' theorem
can be linked to their validity. In the manner of physicists,
it is perhaps not out of place to view these implications as
positive evidence for the structure theory in question.

In relation to the `anabelian' philosophy, we will explain in the
last section how these implications  can be viewed
as a working substitute for the desired implication
`section conjecture $\Rightarrow$ Faltings' theorem.'

\section{The De Rham unipotent Albanese map}

To start out, we let $L$ be any
field of characteristic zero and $S$ be  a  scheme over $L$.
Let $f:X\ra S$ be a smooth scheme over $S$. We denote
by $\Un_n(X)$ the category of unipotent vector bundles with
flat connection having index of unipotency $\leq n$.
That is, the objects are $(\cV,\nabla_{\cV})$,  vector bundles $\cV$
equipped with flat connections
$$\nabla_{\cV}:\cV \ra  \Omega_{X/S}\otimes \cV $$
that admit a
filtration
$$\cV=\cV_{n} \supset \cV_{n-1}\supset \cdots \supset \cV_1 \supset \cV_0=0$$
by sub-bundles stabilized by the connection, such that
$$(\cV_{i+1}/\cV_i, \nabla)\simeq f^*(W_i, \nabla_i),$$
for some bundles with connection $(W_i,\nabla_i)$ on $S$.
The morphisms are maps of sheaves preserving the connection.
Obviously, $\Un_n(X)$ is included in $\Un_m(X)$ as a full subcategory
if $m\geq n$, and we denote by $\Un(X)$ the corresponding union.
Let $b \in X(S)$ be a rational point. It determines
a fiber functor $e_b: \Un(X) \ra \Vect_S$ to the
category of vector bundles  on $S$. $\Un(X)$ forms a Tannakian
category, and we denote by $<\Un_n(X)>$ the Tannakian
subcategory of $\Un(X)$ generated by
$\Un_n(X)$.

Now let $X$ be defined over $L$.
 Given any $L$-scheme $S$ we follow the standard notation of
$X_S$ for the base-change of $X$ to $S$. The points $b$ and $x$ in $X(L)$
then determine fiber functors
$e_b(S), e_x(S): \Un(X_S)\ra \Vect_S$.
We will  use the notation $e_b^n(S)$, etc. to denote the
restriction of the fiber functors to $\Un_n (X_S)$.
The notation $<e_b^n(S)>$ will be used to denote
the restriction to the category $<\Un_n(X)>$,
The {\em De Rham fundamental group} \cite{delignefg}
$\pidr(X,b)$ is the pro-unipotent pro-algebraic group over $L$
that represents the functor on $L$-schemes
$S$:
$$S \mapsto \Aut^{\otimes} (e_b(S))$$
and the path space
$\pidr(X;x,b)$ represents
$$S \mapsto \Isom^{\otimes}(e_b(S), e_x(S))$$
We recall that the isomorphisms in the definition are
required to respect the tensor-product structure, i.e.,
if $g \in  \Isom(e_b(S), e_x(S))$ and
 $V,W$ are objects in $\Un(X_S)$, then
$g(v\otimes w)=g(v)\otimes g(w)$ for $v\in V_b$, $w\in W_b$.
It will be convenient to consider
also the $L$-module $\Hom(e_b, e_x)$ and
the $L$-algebra $\End(e_b)$.

We will fix a base point $b$ and denote the fundamental group
by $U^{DR}$ and the path space by $P^{DR}(x)$. All these constructions
are compatible with base-change.
Given three points $x,y, z\in X$, there is a `composition of paths'
map
$$\pidr(X;z,y)\times \pidr (X;y,x) \ra \pidr(X;z,x)$$
that induces an isomorphism
$$\pidr(X;z,y) \simeq \pidr(X;z,x)$$
whenever one picks a point $p\in \pidr(X;y,x)$ (if it exists). There is an
obvious compatibility when one composes three paths in different orders. In particular,
$P^{DR}(x)$ naturally has the structure of a right torsor
for $U^{DR}$. We denote by $A^{DR}$ and $\cP^{DR}(x)$
the coordinate rings of
$U^{DR}$ and $P^{DR}(x)$, respectively.
As described in \cite{delignefg}, section 10, the $P^{DR}(x)$ fit together to
form the {\em canonical torsor}
$$P^{DR}\ra X$$ which is a right torsor for $X \times_{L}U^{DR} $
and has the property that the fiber over a point $x\in X(L)$ is
exactly the previous $P^{DR}(x)$. We denote by $\cP^{DR}$
the sheaf of algebras over $X$ corresponding to
the coordinate ring of $P^{DR}$.
In order to describe this coordinate
ring, it will be convenient to describe
a universal pro-unipotent pro-bundle associated to
the canonical torsor.

The equivalence (loc. cit.)
from the category of  representations of $U^{DR}$ to unipotent
connections on $X$ can be described as
$$V \mapsto \cV:=(P^{DR}\times V)/U^{DR}$$
Since $P^{DR}(b)$ is canonically isomorphic to
$U^{DR}$, $\cV_b$ is canonically isomorphic to $V$.
Now let $E$ be the universal enveloping
algebra of $\mbox{Lie}U^{DR}$. Then $E$ has the structure of
a co-commutative Hopf algebra and $U^{DR}$
is realized as the group-like elements in $E$. That is,
if $\Delta:E \ra E \hat{\otimes} E$
denotes the co-multiplication of $E$, then
$U^{DR}$ is canonically isomorphic to
$g\in E$ such that $\Delta(g)=g\otimes g$ (of course,
as we vary over points in $L$-algebras).
Let $U^{DR}$ act on $E$ on the left, turning
$E$ into a pro-representation of $U^{DR}$: if $I\subset E$
denotes the kernel of the co-unit of $E$,
then $E$ is considered as the projective system
of the finite-dimensional representations
$E[n]:=E/I^n$.
Define the universal pro-unipotent pro-bundle with
connection on $X$ as
$$\cE:=(P^{DR}\times E)/U^{DR}$$
which is thus given by the projective system
$$\cE[n]:=(P^{DR}\times E[n])/U^{DR}.$$
Then $\cE$ is characterized by the following
universal property:

If $\cV$ is a $\uvb$
and $v\in \cV_b$, then
there exists a unique map
$$\phi_v:\cE \ra \cV$$ such that
$1\in \cE_b \ra v\in \cV_b$.

To see this note that $\cV$ is associated to
a representation $V$ and that a map
$\cE \ra \cV$ is determined by a map
$E=\cE_b \ra \cV_b \ra V$ of representations.
But then, if $1\mapsto v$, then
$f\in E$ must map to $ fv\in V$, so that
the map is completely determined by the image of 1.
To reformulate, there is a natural isomorphism
$$\Hom(\cE, \cV) \simeq \cV_b.$$

By the definition of group-like elements,
the co-multiplication $\Delta$
is a map of $U^{DR}$ representations.
Therefore, there is
a map of connections
$$\Delta: \cE \ra \cE\otimes \cE$$
which turns $\cE$ into a sheaf of
co-commutative co-algebras.
This map can be also characterized as the unique map
$\cE \ra \cE \otimes \cE$ that takes
$1\in \cE_b$ to $1\otimes 1\in (cE\otimes \cE)_b$.
We note that there is a
map
$$\cE_x \ra \Hom (e_b,e_x)$$
defined as follows. Suppose $\cV$ is a $\uvb$,
$v\in \cV_b$ and $f\in \cE_x$. Then
$f\cdot v:=(\phi_v)_x (f)\in V_x$.
It is straightforward to check that
this map is linear in $v$ and functorial in $\cV$.
One also checks that all functorial
homomorphisms $h: \cV_b\ra \cV_x$ arise in this
way: Given such an $h$, consider
its value $f\in \cE_x$ on $1\in \cE_b$. (There is an obvious
way to evaluate $h$ on a pro-$\uvb$. Here as in other places,
we are being somewhat sloppy with this passage.)
Now given any other $v\in \cV$, we have
$\phi_v: \cE \ra \cV$ described above.
By functoriality, we must have a commutative diagram
$$\begin{array}{ccc}
\cE_b & \stackrel{h}{\ra}& \cE_x\\
\scriptsize{\phi_v} \da & & \da \scriptsize{\phi_v} \\
\cV_b & \stackrel{h}{\ra} & \cV_x
\end{array}$$
Hence, we have
$$h(v)=h(\phi_v(1))=\phi_v(h(1))=\phi_v(f)=f\cdot v,$$
proving that the value of $h$ is merely the action of $f$.
Exactly the same argument with $e^n_b, e^n_x$
in place of $e_b, e_x$ shows that
$\cE[n]_x$ represents $\Hom (e^n_b,e^n_x)$.

Now define
$$\cP:=\cE^*,$$  which therefore is
an ind-$\uvb$. The co-multiplication on $\cE$
dualizes to endow $\cP$ with the structure of
a sheaf of commutative algebras over $X$.
Since $P^{DR}(x)$ consists of the
tensor compatible elements in
$\Hom (e_b,e_x)=\cE_x$, we see that this
corresponds to the group-like elements
in $\cE_x$ with respect to the co-multiplication,
and hence, the algebra homomorphisms
$\cP (x) \ra L$. Although we've carried
out this argument pointwise, it applies
uniformly to the whole family as follows:
An argument identical to the pointwise one gives us
maps
$$\cE \ra \underline{\Hom}(\cV_b\otimes \O_X, \cV)$$
of sheaves, that together induce an isomorphism of sheaves
$$\cE \simeq \underline{\Hom}(e_b\otimes \O_X, Id).$$
By definition, $P^{DR}$ represents
$$\underline{\Isom}^{\otimes}(e_b\otimes \O_X, Id).$$
We conclude therefore that $P^{DR}=\Spec_X(\cP)$ and
$\cP=\cP^{DR}$.

Given any pro-algebraic group $G$, we denote by
$Z^nG$ its descending central series normalized by the indexing
$Z^1G=G, Z^{n+1}G=[G,Z^nG]$. Corresponding to this, we have
the quotient groups,
$G_n:=G/Z^{n}G$. When the reference to the
group is clear from the context, we will often omit it from the
notation
and write
$Z^n$ for $Z^nG$. A similar convention will apply to the
various other filtrations occurring the paper. Another convenient
convention we may as well mention here is that we will
often omit the connection from the notation when
referring to a bundle with connection.
That is, $(\cV, \nabla_{\cV})$ will often be denoted simply by $\cV$.

In the case of $U^{DR}$,  we define
$P^{DR}_n(x)$ to be the $U^{DR}_n$ torsor obtained by push-out
$$P^{DR}_n(x)=(P^{DR}(x)\times U^{DR}_n)/U^{DR}$$
Of course it turns out that $P^{DR}_n(x)$ represents the functor
$\Isom^{\otimes} (<e^n_b>,<e^n_x>)$. The pushout
construction can be applied uniformly to the canonical torsor
$P^{DR}$ to get $P^{DR}_n$ which is a
$U^{DR}_n$-torsor over $X$.

If we define $\cP^{DR}[n]:=\cE[n]^*$, it gives a
filtration
$$\cP^{DR}[0]=\O_X\subset \cP^{DR}[1]\subset \cP^{DR}[2]\subset \cdots$$
of $\cP^{DR}$ by finite-rank sub-bundles
that we refer to as the {\em \EM
filtration}.
Observe that the multiplication
$$\cE\otimes E\ra \cE$$
is defined by
taking $f\in E$ and associating to it the unique
homomorphism $\phi_f:\cE \ra \cE$.
In particular, when we compose with the projection
to $\cE[n]$,  it factors through to a map
$$\cE[n]\otimes E[n] \ra \cE[n].$$
Hence, the torsor
map
$$\cP^{DR}\ra \cP^{DR} \otimes A^{DR}$$
carries $\cP^{DR}[n] $ to
$\cP^{DR}[n] \otimes A^{DR}[n]$. Therefore,
if $P^{DR}$ is trivialized by a point $p\in P^{DR}(S)$
in some $X$-scheme $S$, then
the induced isomorphism
$$i_p: \cP^{DR}\otimes \O_S \simeq A^{DR}\otimes \O_S $$
is compatible with the \EM
filtration.
We should also note that if
we embed $L$ into $\C$ and view the whole situation
over the complex numbers, then the \EM
filtration
is the one induced by the length of iterated
integrals \cite{hain2}. In particular, for any
complex point $x\in X(\C)$, the map obtained by base-change
$\cP^{DR}[n]_x\ra \cP^{DR}_x$ is injective.
Thus, the inclusion $\cP^{DR}[n]_x\ra \cP^{DR}_x$
is universally injective. (Proof: Given an $X$-scheme
$f:Z\ra X$, any point $z\in Z$ lies over
$f(z)\in X$. So injectivity of
$$f^*\cP^{DR}[n]_z\ra f^*\cP^{DR}_z$$ reduces
to that of $$\cP^{DR}[n]_{f(z)}\ra \cP^{DR}_{f(z)},$$
which then can be handled through a complex embedding.)
 Therefore, by the constancy of dimension,
 we see that the trivializing map
$i_p$ just mentioned
must be strictly compatible with the \EM
filtration.

Now we let $L=F_v$, the completion of a number field $F$
at a non-archimedean place $v$ and $R_v$ the ring of integers
in $F_v$. Let $X_v$ be a smooth curve over $F_v$.
We will assume that we have a diagram
$$\begin{array}{ccccc}
X_v & \hra &\cX_v &\hra & \cX_v' \\
\da & & \da & & \da \\
\Spec(F_v) & \hra & \Spec (R_v)& = & \Spec(R_v)
\end{array}$$
where $\cX'_v$ is a smooth proper curve over
$R_v$ and $\cX_v$ is the complement in $\cX'_v$
of a smooth divisor $\cD_v$ and $X_v$ is the
generic fiber of $\cX_v$.  We also denote by $X'_v$ (resp. $D_v$)  the
generic fiber of $\cX'_v$ (resp. $\cD_v$).
Letting $k$ denote the residue field of $F_v$, we get
varieties $Y$ and $Y'$ over $k$ as the special fibers of
$\cX_v$ and $\cX_v'$, respectively. Given points
$b, x\in \cX_v(R_v)$, we let $\bar{b}$ and $\bar{x}$ be
their
reduction to $k$-points of $Y$.
Associated to $Y$ and a point $c\in Y (k)$, we have the
{\em crystalline fundamental group} $\picr(Y,c)$ \cite{C-L}
which is a pro-unipotent pro-algebraic
group over $K$, the field of fraction of the ring of
Witt vectors $W$ of $k$.
The definition of this group uses the category $\Un (Y)$ of
over-convergent unipotent isocrystals \cite{berthelot}
on $Y$
which, upon base-change to $F_v$, can  be interpreted as vector bundles
$(\cV,\nabla)$ with connections on $X_v$ satisfying the unipotence
condition of the first paragraph, a convergence condition
on each residue disk $]c[$ of point $c\in Y$, and an over-convergence
condition near the points of $D_v$. For each point $c\in Y (k)$,
we get the fiber functor $e_c:\Un (Y) \ra \Vect_{K}$ which, again
upon base-change to $F_v$, associates
to a pair
$(\cV,\nabla)$ the set of flat sections
$\cV(]c[)^{\nabla=0}$   on the residue disk $]c[$ (\cite{besser}, p. 26).
$\picr(Y,c)$ then represents the  group of tensor automorphisms
this fiber functor. Similarly, using two points and the
functor
$\Isom^{\otimes}(e_c, e_y)$, we get the $\picr(Y,c)$-torsor $\picr(Y;y,c)$
of crystalline
paths  from $c$ to $y$. As before, we will  fix a $c$ in the
discussion and use the notation
$U^{cr}=\picr(Y,c)$, $P^{cr}(y)=\picr(Y;y,c)$.
In fact, if we let $c=\bar{b}$ and $y=\bar{x}$, then we
have canonical isomorphisms (\cite{C-L}, prop. 2.4.1)
$$U^{DR}\simeq U^{cr}\otimes_KF_v$$
$$P^{DR}(x) \simeq P^{cr}(\bar{x})\otimes_KF_v$$
which are compatible with the torsor structures and, more generally,
composition of paths. The key point is that all algebraic unipotent bundles
satisfy the  (over-)convergence condition.
Because  $P^{cr}(y)$ is functorial,
it carries a pro-algebraic automorphism
$$\phi :P^{cr}(y)\simeq P^{cr}(y)$$
induced by the $q=|k|$-power map on $\O_Y$.
The various Frobenius maps are also compatible with the
torsor structures. The comparison isomorphism then endows
$P^{DR}$ with a map that we will again denote by $\phi$.
Given a point  $x \in X_v(R_v)$ that lies in the
residue disk $]\bar{b}[$, we have
$$P^{DR}(x) \simeq U^{cr}\otimes_K F_v$$
and it is worth chasing through the definition to see precisely
what the path is on the left-hand side corresponding
to the identity in $U^{cr}$: Given a unipotent bundle
$(\cV,\nabla)$, we construct an isomorphism
$\cV_b \simeq \cV_x$ by finding the unique flat section $s$ over $]\bar{b}[$
with initial value at $b$ equal to $v_b$. Then the image of
$v_b$ under the isomorphism (i.e., the path) is
$s(x)$, the value of this flat section at $x$.
Since $\phi$ also respects the subcategories
$\Un_n(Y)$ consisting of the overconvergent
isocrystals with index of unipotency $\leq n$,
it  preserves the Eilenberg-Maclane
filtration
on the $\cP^{DR}(x)$.

According to \cite{wojtkowiak}, theorem E, $\cP^{DR}$  possesses
a Hodge filtration. This is a filtration
$$\cP^{DR}=F^0\cP^{DR} \supset \cdots \supset F^n\cP^{DR}
\supset F^{n+1} \cP^{DR} \supset \cdots$$
by sub-$\O_X$-modules that is compatible with
the multiplicative structure. In particular,
the $ F^i$ are ideals. This induces a filtration
on the fibers $\cP^{DR}(x)$ and in particular
on the coordinate ring $A^{DR}$ of $U^{DR}$.  The Hodge filtration
is compatible with the torsor
structure
in that  if the tensor
product $\cP^{DR}\otimes_L A^{DR}$ is endowed with the
tensor product filtration
$$F^n( \cP^{DR}\otimes_L A^{DR})=\Sigma_{i+j=n}  (F^i\cP^{DR})\otimes_L (F^jA^{DR})$$  then the torsor map
$$\cP^{DR} \ra \cP^{DR}\otimes_L A^{DR}$$
is compatible with the filtration.
As in the case of the \EM filtration,
if we utilize their descriptions over $\C$ \cite{H-Z}, we see immediately
that the inclusions $F^i\cP^{DR} \hra \cP^{DR}$ are also universally injective,
and that if we restrict the Hodge filtration to the
terms in the \EM filtration, then the rank of $F^i\cP^{DR}[n]$ is equal to the
dimension of $F^iA^{DR}[n]$.
Since the $F^i$ are ideals, there are corresponding
sub-$X$-schemes
$$F^iP^{DR},$$ where $F^{-i+1}\cP^{DR} $ is the defining ideal for
$F^iP^{DR}$. In fact, $F^0P^{DR}$ is
a $F^0U^{DR}$-torsor.
To see the action, note that compatibility
implies that the torsor map takes
$F^1\cP^{DR}$, the defining ideal for $F^0P^{DR},$
 to $$F^1\cP^{DR}\otimes A^{DR}+\cP^{DR}\otimes F^1A^{DR},$$
the latter being exactly the defining ideal for $F^0P^{DR}\times F^0U^{DR}$
inside $P^{DR}\times U^{DR}$. Thus, we get the action map
$$F^0P^{DR} \ra F^0P^{DR}\times F^0U^{DR}.$$
On the other hand, if we have a point
$p\in F^0P^{DR}(S)$ in some $X$-scheme $S$.
Then the corresponding isomorphism
$$i_p:\cP^{DR} \otimes \O_S \simeq A_{DR} \otimes \O_S$$
is compatible with the Hodge filtration and hence,
strictly compatible  because of the universal injectivity
mentioned above and the equality of ranks (of course
when restricted to the finite-rank terms of the \EM
filtration).
Therefore, it induces an isomorphism
$$i_p:F^0U^{DR}\simeq F^0P^{DR},$$
proving the torsor property. (The remaining compatibilities are obvious
from the corresponding properties for $P^{DR}$.)
Regarding these trivializations, we note
that if a point $p\in P^{DR}(S)$
induces an isomorphism as above compatible with
the Hodge filtration, then since $p$ is recovered
as $ev_0\circ i_p$, where $ev_0:A_{DR} \otimes \O_S \ra \O_S$
is the origin of $U^{DR}_S$, we see that $p\in F^0P^{DR}$.

We will need some abstract definitions corresponding to the
situation described. Given an  $F_v$-scheme $Z$, by a torsor over $Z$ for $U^{DR}$, we mean
 a right $U^{DR}_Z$
torsor $T=\Spec (\cT)$  over $Z$
endowed with an `Eilenberg-Maclane'
 filtration:
$$\cT[0]\subset \cT[1]\subset \cdots$$
and a `Hodge' filtration
$$\cT=F^0T \supset F^1\cT \supset \cdots$$
 of the coordinate ring.
The Eilenberg-Maclane filtration should consist
of locally-free $\O_Z$-modules of finite rank  equal to $\dim_{F_v}A^{DR}[n]$
with the property that $\cT[n]\hra \cT$ is universally
injective.
Furthermore, each $F^i \cT[n]\hra \cT[n]$ must
be universally injective
and have rank equal to
$\dim F^iA^{DR}[n]$.
Both filtrations are required to be compatible with the torsor structure
in slightly different senses, namely,
that
$\cT[n]$ must be carried to $\cT[n]\otimes A^{DR}[n]$
while
$F^m\cT[n]$ is to be taken to
$\Sigma_{i+j=m} F^i\cT\otimes F^jA^{DR}$.
The torsor $T$ must also carry
 a `Frobenius'
automorphism  of $Z$-schemes  $\phi: T\ra T$
(we will denote all of them
by the same letter) which is required to be compatible with
the torsor structure
in the sense that $\phi(t)\phi (u) =\phi(tu)$ for points $t$ of $T$ and
$u$ of $U$. We require that $\phi$ preserves the
Eilenberg-Maclane filtration, although not necessarily
the Hodge filtration.
We will call the
torsor $T$ {\em admissible}, if it is separately trivializable
for the two Frobenius structure and
the Hodge filtration. That is, we are requiring that
there is a point $p^{cr}\in T(Z)$ which is invariant under the Frobenius,
and also that there is a point $p^{DR}\in F^0T(Z)$.

The following lemma follows from the argument of \cite{besser}, cor. 3.2.
 \begin{lem}
Let $T$ any  torsor for $U^{DR}$ over  an affine scheme $Z$. Then $T$ is uniquely trivializable
with respect to the $\phi$-structure. That is, there is a unique
point $p^{cr}_T \in T(Z)$ which is invariant under $\phi$.
\end{lem}
The key fact is that the map
$U^{DR} \ra U^{DR}$, $x\mapsto \phi(x^{-1}) x$ is surjective.
Obviously, the canonical path is the identity for
the trivial torsor $U^{DR}_Z$.

Combining this lemma with the earlier remarks and the unipotence
of $F^0U^{DR}$,
we get that for any affine $X$-scheme $Z$, $P^{DR}_Z$
is an admissible torsor.

We say two admissible torsors are isomorphic
if there is a torsor isomorphism between them
that is simultaneously
compatible with the Frobenius and Hodge filtration.

Let $T$ be an admissible $U^{DR}$ torsor over an $F_v$-algebra
$L$.
Choose a trivialization
$p^H_T \in F^0T (L)$. There then exists a unique element
$u_T \in U^{DR}$ such that
$p^{cr}_Tu_T=p^H_T$. Clearly, a different choice of
$p^H$ will change the result only by right multiplication
with an element of $F^0U^{DR}$. Thus, we get an element
$[u_T] \in U^{DR}/F^0U^{DR}$
which is independent of the choice of $p^H_T$.
\begin{prop}
The map
$$T\mapsto  [u_T]$$
defines a natural bijection
from the isomorphism classes of admissible torsors
to $U^{DR}/F^0U^{DR}$. That is to say, the
scheme $U^{DR}/F^0U^{DR}$ represents the
functor that assigns to each $F_v$ algebra $L$
the isomorphism classes of admissible torsors on $L$.
\end{prop}

{\em Proof.} The map defined is clearly functorial,
so we need only check bijectivity on points.
Suppose
$[u_T]=[u_S]$. Then there exists a $u^0\in F^0U^{DR}$
such that
$u_S=u_Tu^0$.
The elements $p^{cr}_T$ and $p^{cr}_S$ already determine
$\phi$-compatible isomorphisms
$$f: S \simeq U^{DR} \simeq T; \ \ \ p^{cr}_Su \mapsto u \mapsto p_T^{cr}u$$
It suffices to check that this isomorphism is compatible with
the Hodge filtration. But we already know
that the
map
$$h: p^H_Sg \mapsto g \mapsto p^H_Tg$$
is compatible with the Hodge filtration,
and writing the $\phi$-compatible map $f$
with respect to the $p^H$'s,
we get
$$p^{H}_Sg=p^{cr}_Su_Sg\mapsto p^{cr}_Tu_Sg=p^Hu_T^{-1}u_Sg=p^Hu^0g.$$
That is, it corresponds merely to a different choice
$p^Hu^0$ of the  trivialization. So we need to check
that $p^Hu^0\in F^0T$, i.e., that it is a Hodge trivialization.
Its values on the
coordinate ring of $T$ is expressed through the
composition
$$\cT \ra \cT \otimes A^{DR} \stackrel{p^H\times u^0}{\ra} L$$
But then, since the action takes
$F^1\cT$ to $$F^1\cT\otimes A^{DR}+\cT\otimes F^1A^{DR},$$
$F^1\cT$ is clearly killed by this evaluation,
i.e., $p^Hu^0\in F^0T$.

Thus the map of functors described is injective.
To see that it is surjective, note that
for any $Z$, we have $$U^{DR}/F^0U^{DR}(Z)=U^{DR}(Z)/F^0U^{DR}(Z).$$
Now, given $g\in U^{DR}$,
we consider $U^{DR}$ with the same
Hodge filtration, but with the automorphism twisted
to $\phi^g (h)=g^{-1}\phi(g h)$.
This automorphism is compatible with right multiplication, that is,
the right torsor structure on $U^{DR}$:
$$\phi^g (hk)=g^{-1}\phi(g h k)=g^{-1}\phi(g h)\phi(k)=\phi^g(h)\phi(k)$$
Thus, we end up with an admissible torsor $T_g$.
The $\phi^g$-fixed element here is clearly $g^{-1}$,
and thus, $u_{T_g}=g$, as desired.
$\Box$

Already we have the language necessary to define the De Rham unipotent
Albanese map:
$$j^{DR}\cX_v(R_v) \ra U^{DR}/F^0U^{DR} $$
by
$$j^{DR}(x):=[P^{DR}(x)].$$
The finite-level versions
$j^{DR}_n$ are defined by composing with the natural projections
$$U^{DR}/F^0U^{DR} \ra U^{DR}_n/F^0U_n^{DR}.$$
There is a parallel discussion of admissible torsors for
$U^{DR}_n$ out of which one can extract the interpretation of
$$j^{DR}_n:\cX_v(R_v) \ra U^{DR}_n/F^0U_n^{DR}$$
as
$$x\mapsto [P^{DR}_n(x)]$$

In order to describe the map, it will be convenient to have a rather explicit
construction of $U^{DR}$ and $P^{DR}$. We will carry this out assuming $X_v$ is
affine, which is the only case we will need.
The construction  depends on the choice $\a_1, \a_2, \cdots, \a_m$
of global algebraic differential forms representing a basis
of $H^1_{DR}(X_v)$. So here, we have $m=2g+s-1$, where
$s$ is the order of $X'_v\setminus X_v$.
Corresponding to this choice, there is a
 free non-commutative algebra
$$F_v<A_1, \ldots, A_{m}>$$
generated by the symbols
$A_1, A_2, \ldots, A_{m}$.  Thus, it is the group algebra of
the free group on $m$ generators. Let
$I$ be its augmentation ideal.
The algebra
$F_v<A_1, \ldots, A_{m}>$ has a natural comultiplication
map $\Delta$ with values
$\Delta (A_i)=A_i \otimes 1+ 1\otimes A_i$.
Now we consider the completion
 $$E:=\invlim F_v<A_1, \ldots, A_m>/I^n$$
 $\Delta$  extends naturally to a comultiplication
$E \ra E\hat{\otimes} E$.
We can also consider the quotients
$$E[n]:=E/I^{n+1}A=F_v<A_1, \ldots, A_{m}>/I^{n+1}$$
that carry the induced comultiplication.
Thus, $E$ and $E_n$ have the structure of non-commutative,
co-commutative Hopf algebras.

Now let $\cE$ be the pro-unipotent pro-vector bundle
$E\otimes \O_{X_v}$ with the connection $\nabla$
determined by
$$\nabla f=-\Sigma_i A_if\a_i$$
for  constant sections $f\in E$.
We also have the finite-level quotients
$\cE[n]:=E[n]\otimes \O_{X_v}$.
This construction ends up as the `universal' pro-unipotent
pro-bundle with connection in the sense of the beginning paragraphs,
justifying the conflation of notation.
To see this, we need a lemma.

\begin{lem}
Let $(V,\nabla)$ be a unipotent bundle with flat connection on
$X_v$ of rank $r$. Then
there exist strictly upper-triangular matrices $N_i$
 such that
$$(V,\nabla) \simeq (\O_X^r, d+\Sigma_i \a_i N_i)$$
\end{lem}
{\em Proof of lemma.}
Since $X_v$ is affine and $(V,\nabla)$ is unipotent, by choosing
vector bundle splittings of the filtration, there exists
a trivialization
$V\simeq \O_X^r$ such that
$\nabla$  takes the form
$d+\om$ for some strictly upper-triangular
$n\times n$ matrix of 1-forms $\om$. It will be convenient to write
$\om=\Sigma \om_{ij}E_{ij}$
where $E_{ij}$ is the elementary $n\times n$ matrix with
a 1 in the $(i,j)$-entry and zero elsewhere, and $\om_{ij}=0$
unless $j>i$.
Recall that a gauge transformation $G$ will change
the connection matrix by
$$\om \mapsto G^{-1}\om G +G^{-1}dG$$
For a gauge transformation of the form
$G=I-aE_{ij}$ ($j\neq i$), we have $G^{-1}=I-aE_{ij}$
while $G^{-1}dG= daE_{ij}$. We wish to perform a series of gauge
transformations so that each entry $\om_{ij}$
is replaced by a linear combination of the $\a_i$'s and $\b_j$'s.
Now, any single one-form can be written as such a combination
after adding an exact form.
 We will show how to change $\om_{ij}$ by induction on
$j-i$. Assume we are done for $j-i<c$. So let $i,j$
satisfy $j-i=c$.
First find $a$ such $\om_{ij}+da$ is a linear combination
of the $\a_i$ and $\b_j$.
Now consider the gauge transformation
$G=I-aE_{ij}$. We compute
$$G^{-1}\om=(I+aE_{ij})(\Sigma \om_{kl}E_{kl}) \\
=\om+\Sigma_l \om_{jl}aE_{il}=:\om'$$
Note the only non-zero $l$ occurring in the sum are strictly
bigger than $j$. Thus, $\om'_{\mu\nu}\neq \om_{\mu\nu}$
only occurs for $\mu=i$ and $\nu >j$, that is, for
$\nu-\mu>j-i$.
On the other side,
$$\om'G=(\Sigma \om'_{kl}E_{kl})(I-aE_{ij})\\
=\om'-\Sigma \om'_{ki}E_{kj}$$
and all the $k$ appearing in the sum are $<i$.
Therefore, we see that all the entries
$\om^G_{kl}$ of  $\om^G:=G^{-1}\om G$ are
equal to $\om_{kl}$ for $l-k\leq c$.
On the other hand, $dGG^{-1}=daE_{ij}$.
Performing such a gauge transformation for each
$i,j$ such that $j-i=c$ clearly achieves what we want.
$\Box$.

Now we can formulate the universal property of $\cE$:
\begin{lem}
Given any object  $(V,\nabla)$ in $\Un(X_v)$ together with a
an element $v \in V_b$ (the fiber at $b$), there exists a unique morphism
$\phi: \cE \ra V$ in $\Un(X_v)$ such that $ 1\in \cE_b\mapsto v$.
That is, $\cE$ is the same as the bundle with the
same notation from the beginning of the section.
\end{lem}
{\em Proof.}
First we will show uniqueness. For this, we can choose an embedding
of $F_v$ into the complex numbers and assume that everything is
defined over $\C$ (and for convenience, we will omit the base-change to $\C$ from the notation).
By uniqueness of solutions to differential
equations with given initial condition, we see that
any map of bundles $\phi: \cE \ra V$ is completely determined by
its value at $b$. (Actually, for this, we needn't go to $\C$.
It suffices to use uniqueness of formal solutions.)
So we need to check that the value of the map is determined
at $b$ by the given conditions. Let  $\pi$
be the topological fundamental group of
$X_v(\C)$ and let $\C[\pi]$ be the complex group algebra of $\pi$.
The holonomy transformations give an action of
$\C[\pi]$ on the fibers $\cE_b$ and $V_b$. Since $\phi$ respects
the connections,
$\phi_b$ is equivariant for this action.
But over $\C$, one can compute the holonomy using iterated
integrals. Given  a word $w=A_{i_1}A_{i_2}\cdots A_{i_n}$ in the $A_i$ and a piece-wise
smooth path
$\g$,
denote by $\int_{\g} \a_w$, the iterated integral
$$\int_{\g} \a_{i_1}\a_{\i_2} \cdots \a_{i_n}.$$
We are using the normalization whereby if
$\g^*(\om_i)=f_i(t)dt$, then
$$\int_{\g} \om_1\om_2 \cdots \om_n=\int_0^1f_1(t_1)
\int_0^{t_1}f_2(t_2)\int_0^{t_2}f_3(t_3)\int \cdots \int_0^{t_{n-1}}f_n(t_n)dt_n dt_{n-1}\cdots dt_1.$$
On a contractible open set $U$, if we pick a base point
$c$ and define the function $G_c(z)$ with values in $E$ on $U$
by
$$G_c(z)=\Sigma_w \int_c^z \a_w[w],$$
where the iterated integral occurs along any path from $c$ to
$z$, then
$dG_c=\Sigma_i\a_i A_iG_c$
so that $G_c(z)$ can be viewed as a flat section of $\cE$.
In particular, for the holonomy around a loop $\g$ based at $b$,
$$\g\cdot f=\Sigma_w \int_{\g} \a_w[w]f$$
 for $f\in E=\cE_b$.
Thus, we must have
$$\phi (\Sigma_w \int_{\g} \a_w[w]1)=\g\phi(1)=\g v$$
Now consider the map
$\C[\pi]\ra E$ given by the formula
$$\g \mapsto G_{\g}:=\Sigma_w \int_{\g} \a_w[w]$$
(We will use the natural notation
$\int_{\g}\a_w:=\Sigma_ic_i\int_{\g_i} \a_w$ if
$\g=\Sigma_ic_i\g_i$, $c_i\in \C, \g_i\in \pi$.)
To conclude the proof of uniqueness, we just need
to see that the composition to each of the finite-dimensional
quotients $E_n$ is surjective. But if not, we would have a non-trivial
linear relation between the
$\int_{\cdot} \a_w$ regarded as functions on $\C[\pi]$.
We show by induction on the length of $w$ that they are in fact
linearly independent.
For this, let $a_1, \ldots, a_m$ be elements of $\C[\pi]$
whose classes in homology form a dual basis to
the $\a_i$. Assume
$$\Sigma_{|w| \leq n} c_w \int_{\g} \a_w=0$$
for all $\g \in \pi$. Write this as
$$\Sigma_{|w|\leq n-1} ( \Sigma_ic_{A_iw} \int_{\g} \a_{A_iw})=0$$
and apply it to $a_j\g$. Then we have
$$\Sigma_{|w|\leq n-1} ( \Sigma_ic_{A_iw} \int_{a_j\g} \a_i\a_w)=0$$
The co-product formula for iterated integrals \cite{chen} can then be applied togive
$$\Sigma_{|w|\leq n-1} [ \Sigma_ic_{A_iw} (\mbox{deg}(a_j)\int_{\g} \a_i\a_w+\int_{a_j}\a_i \int_{\g}\a_w+\cdots)]=0,$$
where
$\mbox{deg}(\Sigma c_{\g} [\g])=\Sigma c_{\g}$ on
elements of $\C[\pi]$.
From this, we subtract ($\mbox{deg}(a_j)$ times)
the original relation to get
$$\Sigma_{|w|=n-1}c_{A_jw}\int_{\g}\a_w+\Sigma_{|w|\leq n-2}b_w\int_{\g}\a_w=0$$
for some constants $b_w$. Since this holds for all $\g$, by induction, we have
$c_{A_jw}=0$ for $|w|=n-1$, and hence $c_w=0$ for $|w|=n$ when the argument is applied to each $A_j$.
Then again by induction, all the $c_w$ must vanish.

The existence part is
easy. We may assume that $(V,\nabla)$ is of the form
$(\O_{X_v}^r, d-\Sigma_i\a_iN_i )$ as in the previous lemma, except we have changed $N_i$ to $-N_i$.
Now given the vector $v\in V_b$, define the map
$\cE \ra V$ that sends a section
$\Sigma f_w[w]$ of $\cE$ to $\Sigma_w f_w N_wv$ where the notation is
that if a word has the form
$w=A_{i_1}A_{i_2}\cdots A_{i_n}$, then $N_w$ is the matrix
$N_{i_1}N_{i_2}\cdots N_{i_n}$. This map is obviously
compatible with the connection and takes the constant section 1
to the constant section $v$ (in particular, at the point $b$).
$\Box$

As mentioned already,
another way of phrasing the lemma is to say that there is a canonical
isomorphism:
$$\Hom(\cE,V)\simeq V_b$$
In particular, $\Hom(\cE,\cE)\simeq \cE_b=E$
and it is easily checked that this isomorphism is
compatible with the algebra structure.

Recall the evaluation fiber functor $e_b:\Un(X_v)\ra \Vect_{F_v}$.
Since $\cE$ is the universal pro-unipotent
vector bundle with connection, $E=\cE_b$  represents
$\End(e_b)$ and
$E=\cE_x$ represents $\Hom (e_b,e_x)$.
It is equally easy to see that
$E[n]$ represents $\End (e^n_b)$ and $\Hom (e^n_b,e^n_x)$.

If we examine
the co-multiplication
$\Delta :\cE\ra \cE\otimes \cE$
associated to the $\cE$ as a universal object,
that is, the unique map that sends
1 to $1\otimes 1$, by the definition of the connection on $\cE$,
it takes $A_i$ to $A_i\otimes 1+1\otimes A_i$. That is,
the co-multiplication is compatible with the
one defined formally on $E=\cE_b$.
So we also have that
 $\Aut^{\otimes}(e_b)$ and $\Isom^{\otimes}(e_b, e_x)$ are represented
by the group-like elements in $\cE_b=E=\cE_x$
and that $\cE^*$, the dual bundle of $\cE$, with its structure
of a commutative algebra is nothing but $\cP^{DR}$.

Now $\cP^{DR}$ carries a Hodge filtration by sub-bundles,
but its explicit description will be unnecessary for our purposes.
All we will need is that $F^0P$ is a torsor for the
group $F^0U^{DR}$, and hence, can be trivialized over $X_v$
in the case it is affine.

{\em Proof of theorem 1.}

First note that if $X_v$ happens to be compact,
we can delete an $R_v$ point $x_0$ away from the residue disk of
$b$ to get $\cX^0_v=\cX_v\setminus x_0$, still equipped with an $R_v$ point $b$. Since restriction
embeds the unipotent bundles with flat connection on
$X_v$ as a full subcategory  of the unipotent
bundles on $X^0_v$ (see, for example,
\cite{besser}, cor. 2.15), we get a surjection
$\pidr (X^0_v,b) \ra \pidr (X_v,b)$. This map is compatible
with the Frobenius because the comparison isomorphism
with $\pi_1^{cr}$ is functorial \cite{C-L}.
Functoriality of the Hodge filtration can be deduced
either from an examinination of Wojkowiak's construction or the usual
functoriality over $\C$ \cite{hain2}.  Of course, an
identical argument
can be applied to the path spaces,
giving us a commutative diagram
$$\begin{array}{ccc}
\cX^0(R_v) & \ra & \pidr (X^0_v,b)/F^0\\
\da & & \da \\
\cX(R_v) & \ra & \pidr (X_v,b)/F^0\\
\end{array}$$
where the right vertical map is surjective.
Therefore, it suffices to prove the theorem for
an affine curve.

In fact, we will show that the image of the residue disk
$]\bb[$ is already Zariski dense.
Now, as described above,
there exists a canonical $U^{DR}$ torsor $P^{DR}$
on $X_v$ with the property that the fiber over $x$ is exactly
$P^{DR}(x)$ which is, in fact, admissible. This is because
of the unique Frobenius invariant element
and the fact that $F^0P^{DR}$ is a $F^0U^{DR}$-torsor,
and hence, has a section over an affine scheme.
  Choose such a trivialization
$p^H \in F^0P^{DR}(X_v)$. It is important to note here that $F^0P^{DR}$ becomes thereby
{\em algebraically} trivialized.
We can compare this with
any other algebraic trivialization
$g\in P^{DR}(X_v)$. That is, there exists an algebraic map
$\g:X_v \ra U^{DR}$ such that
$g\g=p^H $.  For  any $x\in X_v$, we have a point
$u(x)\in U^{DR}$ such that $p^{cr}(x)u(x)=g(x)$so
that $p^{cr}(x)u(x)\g(x)=p^H(x)$.
That is,
$j^{DR}(x)=[u(x)\g (x)]$. Therefore, it suffice to
show that the image of
$x\mapsto u(x)\g (x)$ is Zariski dense in $U^{DR}(X)$.
 As previously discussed, given $(V,\nabla)$, on
$]\bb[$, the element $p^{cr}(x)$
is obtained from the diagram
$$\begin{array}{ccccc}
 & & V(]\bb[)^{\nabla=0} & & \\
 & \scriptstyle{\simeq}\swarrow  & & \searrow \scriptstyle{\simeq}& \\
V_b & & & & V_x
\end{array}$$
as the inverse of the left arrow followed by the right arrow.
We proceed now to describe the various objects in local coordinates
to obtain the desired result.

 Take $g$ to be the trivialization of $P^{DR}$
 determined by the previously described trivialization of $\cE$.
 That is, since we have identified
 $\Isom^{\otimes}(e_b,e_x)=E$, the element
 in question is $1\in E$.
 For $x \in]\bb[$ $u(x)$ is then
identified with the value of the function
 $u:]\bb[ \ra E$ that is the unique
horizontal section of $\cE$ convergent on $]\bb[$
such that $u(b)=1$. Write
$u=\Sigma_wu_w[w]$, where $w$ runs over the words in $A_i$.
We will show:
\begin{lem}
The $u_w$ are linearly independent over the algebraic
functions on $X_v$.
\end{lem}

It is easy to see that this lemma implies the theorem:

We have $$j^{DR}(x)=\Sigma_w j_w(x)[w]=
u(x)\g(x)=\Sigma [\sigma_{w'w''=w}u_{w'}(x)\g_{w''}(x)][w],$$
where $\g=\Sigma_w\g_w[w]$
If the image of $j^{DR}_n$ were in a Zariski closed subspace,
we would have a linear relation among the $j_w$, by the
description we have given of the coordinate ring of
$U^{DR}$.
Let us see by induction on the length of $w$
that the $j_w$ are linearly independent.
Assume we have a relation
$$\Sigma_{|w|\leq n} c_w j_w=0$$
Then
$$\Sigma_{|w|\leq n} c_w[\Sigma_{w'w''=w}u_{w'}\g_{w''}]=0.$$
This would give us a relation of the form
$$\Sigma_{|w|\leq n} a_w u_w=0$$
where the $a_w$ are algebraic functions. Hence, all $a_w=0$.
But for $|w|=n$, we have
$a_w=c_w\g_0$ and $\g_0=1$  because $\g$ is group-like.
Therefore, we must have $c_w=0$ for all $|w|=n$. The result follows by induction.

So we need only give the

{\em Proof of lemma.} The proof here is entirely similar to
the uniqueness part of Lemma 4.
Even though the $u_w$ are convergent on $]\bb[$ we can regard
them as formal power series and show independence over the
algebraic functions in there. For this, we choose an
embedding $F_v\hra \C$ into the complex numbers,
so that we are in the situation of a connection on
a curve over $\C$. Uniqueness of formal solutions
show that the $u_w$ are just the power series expansions
of analytic solutions near $b$. But the analytic solutions can
be obtained as iterated integrals:
in a contractible neighborhood of $b$, choose a path $c$
from $b$ to $x$. Then for
$w=A_{i_1}A_{i_2} \cdots A_{i_n}$, we have
$$u_w=\int_c \a_{i_1}\a_{i_2}\cdots \a_{i_n}$$
We see therefore, that $u(x)$ can also be continued
as a multi-valued function to all of $X_v$ (which, recall,
we are now regarding as a complex curve). Regard the
$u_w$ as functions on the universal covering
$\tilde{X}_v$ after choosing a point $b'\in \tilde{X}_v$
lying over $b$. Let the topological fundamental group
$\pi$ of $X_v$ act this time  on the space of analytic
functions on $\tilde{X}_v$ by
$\g: f(z) \mapsto f^{\g}(z)=f(z\g)$. This extends to an action of the
complex group algebra $\C[\pi]$. To compute $u_w^{\g}(z)$
for an element $\g \in \G$, we proceed
as follows: Choose a path $c$ from $b'\g$ to $z$. Now choose a lifting
$\tilde{\g}$ of $\g$ to a path
in $\tilde{X_v}$ from $b'$ to $b'\g$. Then
$$u_w(z\g)=\int_{c\tilde{\g}} \a_{i_1} \cdots \a_{i_n} \\
=u_w(z)+\int_{c}\a_{i_1} \cdots \a_{i_{n-1}}
\int_{\tilde{\g}}\a_{i_n}
+\int_{c}\a_{i_1} \cdots \a_{i_{n-2}}\int_{\tilde{\g}}\a_{i_{n-1}}\a_{i_{n-2}}+\cdots
$$
 As before, let $a_i\in \C[\pi]$ be elements whose
image in $H_1(X_v,\C)$ form a basis dual to the $\a_i$.
Then we see that
$$
u_{wA_i}^{a_i}=u_{wA_i}+u_w+\Sigma_{|w'|<|w|}c_{w'}u_{w'}$$
for constants $c_{w'}$.

Now we prove the desired independence by induction on
$|w|$. So assume we are done for $|w|<n$, and suppose we have
a linear relation
$$\Sigma_{|w|\leq n} h_wu_w=0$$
where the $h_w$ are algebraic. In particular, they are
single-valued on $X_v$, and so are acted on trivially
by $\pi$. We write this relation as
$$\Sigma_{|w|= n}h_w u_w+\Sigma_{|w|=n-1}h_w u_w+\Sigma_{|w|\leq n-2}h_w u_w=0$$or
$$\Sigma_i [\Sigma_{|w|= n-1}h_{wA_i} u_{wA_i}+\Sigma_{|w|=n-2}h_{wA_i} u_{wA_i}]+\Sigma_{|w|\leq n-2}h_w u_w=0$$
When we apply $a_j-I$ to this relation, we get
$$\Sigma_{|w|=n-1}h_{wA_j}u_{w}+\Sigma_{|w|\leq n-2}b_wu_w=0$$
for some algebraic functions $b_w$.
By induction, we get $h_{wA_j}=0$ for all $w$. Since this works
for any $j$, we are done.
$\Box$

\section{\'Etale Realizations}

This time we start with a smooth scheme $X$ over a field
$L$ of characteristic zero and put
$\bX=X\otimes_L \bar{L}$. We have on
$\bX$ the category $\Un^{et}_p(\bX)$
of $\Q_p$-unipotent lisse sheaves.
Choosing a point $b\in \bX$ then
determines a fiber functor
$e_b:\Un^{et}_p(X) \ra \Vect_{\Q_p}$
by taking stalks. We define the pro-unipotent
$p$-adic \'etale fundamental group (see, e.g., \cite{wildeshaus})
and path space by
$$\piet(\bX,b)=\Aut^{\otimes}(e_b)$$
and
$$\piet(\bX;x,b)=\Isom^{\otimes} (e_b,e_x),$$
where, as usual, we have represented an obvious functor
on $\Q_p$ algebras. When both $b$ and $x$ are in $X(L)$,
we get an action of $G=\Gal (\bar{L}/L)$.
In particular, $P^{et}(x)=\piet(\bX;x,b)=\Spec(\cP^{et})$
is a right $U^{et}=\piet(\bX,b)=\Spec (A^{et})$ torsor in a manner compatible
with the Galois action. Here as well, we can
consider the subcategory $\Un^{et}_{p,n}(\bX)$
of unipotent $\Q_p$ local systems $V$ having
index of unipotency $\leq n$ in that $V$ is
endowed with a filtration
$$V=V_{n}\supset V_{n-1} \supset \cdots \supset V_1\supset V_0=0$$
such that each successively quotient is isomorphic to
a trivial $\Q_p$-sheaf. The corresponding restrictions
of the fiber functors will again  be denoted
$e^n_b$, etc. We denote by
$<\Un^{et}_{p,n}(\bX)>$ the Tannakian subcategory of
$\Un^{et}_p(\bar{X})$ generated by
$\Un^{et}_{p,n}(\bX)$. Using
the restrictions $<e^n_b>:=e_b|<\Un^{et}_{p,n}(\bX)>$
and $<e^n_x>:=e_x|<\Un^{et}_{p,n}(\bX)>$
 of the fiber functors, we then
get $U_n^{et}$ and $P^{et}_n(x)$.
In the above, we can everywhere consider homomorphisms of functors
rather than tensor isomorphisms, like
$\Hom(e_b,e_x)$,  which of course have
the structure of  $\Q_p$ vector spaces.

It might be useful to recall at this point some rudimentary points about the
structure of these groups and torsors, although it is not necessary
  to be
as explicit as in the De Rham case.

Let $E^{et}$ be the universal enveloping algebra of
$\mbox{Lie }U^{et}$. Then $E^{et}$ again has the structure
of a co-commutative Hopf algebra, and
we define
$E^{et}[n]:=E^{et}/I^n$, where $I\subset E$ is
the kernel of the counit $E^{et}\ra \Q_p$.
As in the De Rham situation, we have the
universal pro-unipotent lisse $\Q_p$-sheaf $\cE^{et}$
that can be constructed, for example, by
twisting the canonical profinite torsor
$\hat{T}:=P^{\wedge}|X\times b$, in the notation of \cite{delignefg}, section 10.17,
with representation of the profinite $\hat{\pi}_1(\bar{X},b)$ obtained
by composing
the natural map
$$\hat{\pi}_1(\bar{X},b) \ra U^{et}$$
with the left multiplication of $U^{et}$ on $E$.
By an argument identical to the De Rham setting,
we deduce the property that $\cE^{et}_x$
represents $\Hom(e_b,e_x)$ and the finite-rank quotients
$\cE^{et}[n]$ corresponding to $E/I^n$ represent
in their turn $\Hom (e^n_b, e^n_x)$.
Once again the universal property gives
us a comultiplication
$$\Delta:\cE^{et} \ra \cE^{et}\otimes \cE^{et}$$
(here and henceforward, a tensor product
without an explicit subscript will refer to tensoring over
$\Q_p$) so that if we define
$\cP^{et}=\Hom (\cE^{et}, \Q_p)$, then
$\cP^{et}$ is an ind-$\Q_p$ lisse sheaf with the
property that $\cP^{et}_x$ is the coordinate ring
of $P^{et}(x)$. There is also a parallel (to the De Rham setting)
discussion
corresponding to $\cP^{et}[n]$, $\cP^{et}_n$,
and $P^{et}_n$, which we will assume without further comment.
We will refer again to the filtration by the
$\cP^{et}[n]$ as the \EM filtration.
When $x\in X(L)$, since the categories
$\Un^{et}_{p.n}(\bar{X})$
are invariant under the Galois action,
so is
$\cP^{et}[n](x)$.
Therefore, $\cP^{et}$ is a
direct limit of the finite-dimensional Galois modules
$\cP^{et}[n](x)$,
allowing much of the theory of finite-dimensional representations,
for example Fontaine's theory of the Dieudonn\'e functor \cite{fontaine},
to apply.

In the subsequent discussion, fix
 a quotient group $G$ through which the action of $\Gal(\bar{L}/L)$
factors.

We recall some notions from \cite{kim}.
Given a $\Q_p$-algebra $R$,
we give it the inductive limit topology obtained from the
natural topology on finite-dimensional $\Q_p$-subspaces.
Then by a torsor for $U^{et}$ over $R$, we mean a $U^{et}_R$
torsor $T=\Spec (\cT)$ in the usual sense,
except that we require that

(1)$\cT=\dirlim_n \cT[n]$ for some locally free sub-$R$-modules
$\cT[n]$ that are local direct summands and
of finite rank equal to $\dim A^{et}[n]$,
 as $n$ runs through the non-negative integers.

(2) $\cT$ is equipped with
a continuous $R$-linear action of $G$  that stabilizes the $\cT[n]$.

(3) The torsor structure
$$\cT \ra \cT \otimes_R (A^{et}\otimes R)$$
is compatible with the G-action and takes
$\cT[n]$ to
$\cT[n]\otimes_R(A^{et}[n]\otimes R)$.

When $L$ is the completion $F_v$ of a number field
$F$ as in the introduction
 and $G=G_v=\Gal(F_v/F)$, we get  a classifying pro-variety
$H^1(G_v, U^{et})$ \cite{kim} defined over $\Q_p$ for  torsors, as well as
a subvariety $H^1_g(G_v,U^{et})$, consisting of those
torsors
that trivialize upon base-change to $B_{DR}$. We warn the reader at this point that
 this subvariety was denoted $H^1_f$ in the reference
\cite{kim}. In this paper, we will use the notation
$H^1_f$ for a certain subvariety of $H^1_g$
defined using $B_{cr}$ in place of $B_{DR}$
(see the following section).
Let us spell out what this $B_{DR}$ condition means.
If $T=\Spec (\cT)$ is a torsor for $U^{et}$ over $R$, we say
it is a {\em De Rham torsor} if there is a $G_v$-equivariant
algebra homomorphism
$$t:\cT \ra R\otimes B_{DR}$$
or equivalently,
$$t_{B_{DR}}: \cT\otimes B_{DR} \ra R\otimes B_{DR}.$$
That is, we are requiring the existence
of a $G_v$-invariant point in $T_{B_{DR}}(R\otimes B_{DR})$.
Then $t$ induces an isomorphism of $B_{DR}$-algebras
$$C_t: \cT\otimes B_{DR} \simeq A^{et}\otimes R\otimes B_{DR}$$
that is $G_v$-equivariant. Furthermore, since
$t_{B_{DR}}$ is a homomorphism of
$B_{DR}$ algebras, it respects the Hodge filtration
induced by that on $B_{DR}$. Hence, so does
$C_t$.
Since
there is  the identity homomorphism
$$e: A^{et}\otimes R \otimes B_{DR}\ra R\otimes B_{DR} ,$$
out of which we recover $t$ as $e\circ C_t$,
we see that the datum of $C_t$  is
equivalent to the De Rham condition.

Taking $G_v$-invariants,
we get an isomorphism
$$C_{t^{G_v}}: D(\cT):=(\cT\otimes B_{DR})^{G_v} \simeq  (A^{et}\otimes B_{DR})^{G_v}
\otimes R
\simeq A^{DR}\otimes R=A^{DR}\otimes_{F_v} (F_v\otimes R)$$
Here, we are using the fact
 that $D(A^{et})\simeq A_{DR}$ in a manner compatible
with the Hodge filtration \cite{olsson}.
As a consequence of the definition, the
coordinate ring $\cT$ becomes a limit of finite-dimensional
De Rham representations in the sense of
\cite{fontaine} and hence, the inclusion
$$D(\cT) \hra \cT\otimes B_{DR}$$
induces an isomorphism
$$D(\cT)\otimes_{F_v} B_{DR}\simeq \cT\otimes B_{DR}.$$
If we  examine what happens to the torsor structure
$\cT \ra \cT \otimes_R(A^{et}\otimes R),$
it extends to
$$\cT\otimes B^{DR} \ra \cT \otimes_R(A^{et}\otimes R) \otimes B^{DR}$$
But
$$
\cT \otimes_R(A^{et}\otimes R) \otimes B^{DR}=(\cT\otimes A^{et})\otimes B_{DR}
=(\cT\otimes B_{DR})\otimes_{B_{DR}} (A^{et}\otimes B_{DR})$$
$$
=((\cT\otimes B_{DR})^{G_v}\otimes_{F_v} B_{DR})\otimes_{B_{DR}}
((A^{et}\otimes B_{DR})^{G_v}\otimes_{F_v} B_{DR})=D(\cT)\otimes_{F_v}
D(A^{et})\otimes B_{DR}
$$
So taking $G_v$-invariants gives us
$$D(\cT) \ra D(\cT) \otimes_{F_v} A^{DR}=
D(\cT)\otimes_{R\otimes F_v} (A^{DR}\otimes_{F_v} (R\otimes F_v)).$$
This defines the structure of a $U^{DR}$-torsor over $R\otimes F_v$
on $D(T)$ which is trivialized by
the point $t^{G_v}: D(\cT) \ra R\otimes F_v$ obtained from
$t$ by restricting to the $G_v$-invariants. In fact,
the map
$$D(\cT)\ra A^{DR}\otimes_{F_v}(R\otimes F_v)$$
obtained from $t^{G_v}$ is seen to be none
other than the earlier isomorphism
$C_{t^{G_v}}$. The  compatibility
of the torsor structure with the Eilenberg-Maclane filtration
follows from the fact that the comparison isomorphism
of \cite{olsson} arises from an equivalence of
categories $$\Un^{et}_{p}(\bar{X})\otimes B_{DR} \simeq \Un^{DR}(X)\otimes_{F_v}B_{DR}$$
that respects the subcategories specified by the index of unipotency
on either side.

The Hodge filtration on $D(\cT)$ is defined to be the one induced by
the Hodge filtration  on $B_{DR}$. Since this is
true also of the one on $A^{DR}=D(A^{et})$
and
$\cT\otimes B_{DR}\simeq A^{et}\otimes R\otimes B_{DR}$
preserves the filtration, so does
$D(\cT) \simeq A^{DR}\otimes_{F_v}R$.
We conclude that the
$U^{DR}$ torsor thus obtained can be trivialized
 together with the Hodge filtration,
in the sense of the previous section.

There is a parallel discussion with $B_{cr}$
in place of $B_{DR}$ that defines a {\em crystalline}
$U^{et}$-torsor $T$ over $R$ to be one that trivializes over
$B_{cr}$ or, equivalently, admits a $G_v$-equivariant
isomorphism of torsors
$$\cT \otimes B_{cr} \simeq A^{et}\otimes B_{cr}.$$
Using the usual  inclusion of $B_{cr}$ into $B_{DR}$,
we see that a crystalline torsor is also De Rham.
It is easy to show that the crystalline condition
defines a subvariety
$$H^1_f(G_v, U^{et})\subset H^1_g(G_v, U^{et})$$
However, since the reference \cite{kim} dealt explicitly
only with $H_g$,
we sketch the modification necessary to represent the
crystalline torsors. For each $n$, consider the functor
on $\Q_p$-algebras
$$H^0(G_v, U^{et,B_{cr}}_n/U^{et}_n): R\mapsto H^0(G_v, U^{et}_n(B_{cr}\otimes R)/
U^{et}_n(R))$$
Then we have
\begin{lem}
$H^0(G_v, U^{et,B_{cr}}_n/U^{et}_n)$
is representable by an affine variety over $\Q_p$.
\end{lem}
The proof is verbatim the same as \cite{kim}, section 1, proposition
3 and will be therefore omitted.
And then, as in the few paragraphs preceding that
proposition, we have an exact sequence
$$0\ra U^{et}_n(R) \ra U^{et}_n(R\otimes B_{cr}) \ra U^{et}_n(R\otimes B_{cr})/U^{et}_n(R)\ra 0$$
from which we get the connecting homomorphism
$$H^0(G_v, U^{et}_n(R\otimes B_{cr})/U^{et}_n(R)) \ra H^1(G_v, U^{et}_n(R))$$
which is functorial in $R$. Thus, it is represented by a map of
varieties, and its image is a subvariety of
$H^1(U^{et}_n(R))$. But by looking at the crystalline
condition as the existence of a $G_v$-invariant point in
$P^{et}(R\otimes B_{cr})$,
we see that the image consists exactly of the crystalline
torsors (exactly as in the proof of
\cite{kim}, section 2, lemma 6). By passing to the limit over $n$, we see that
the set of
isomorphisms classes of crystalline torsors for $U^{et}$
is represented by a pro-algebraic subvariety $H^1_f(G_v,U^{et})$
of $H^1_g(G_v,U^{et})$.

If we denote by $K$ the absolutely unramified subfield of
$F_v$, then we get a functor
$D_{cr}$ from crystalline $U^{et}$-torsors over $R$ to (trivial)
$U^{cr}$-torsors over $R\otimes K$ endowed with a Frobenius endomorphism:
This functor of course sends
$T=\Spec(\cT)$ to $$D_{cr}(T):=\Spec ((\cT\otimes B_{cr})^{G_v}).$$
A discussion entirely parallel  to the $B_{DR}$ case
shows that $D_{cr}(T)$ is endowed with a Frobenius endomorphism
and is a torsor for $U^{cr}=D_{cr}(U^{et})$ in a manner compatible with
this extra structure. As always, this torsor
is canonically trivial.
We have the comparison
$$D(T)\simeq D_{cr}(T)\otimes_K F_v.$$
 Thereby, for crystalline torsors $T$, $D(T)$ is endowed with
both a Hodge filtration and a Frobenius.

To summarize,
 the functor $T\mapsto D(T)$ defines a map from
the isomorphism classes of crystalline $U^{et}$-torsors over
$R$
to admissible $U^{DR}$ torsors over $R\otimes F_v$.
It is therefore represented by
a pro-algebraic map

$$D: H^1_f(G_v, U^{et}) \ra Res^{F_v}_{\Q_p} (U^{DR}/F^0).$$

Going back from abstract considerations to
the study of rational points, each $P^{et}(x)$ defines a class in
$H^1_f(G_v,U^{et})$ \cite{olsson}. We are thus given  an \'etale unipotent
Albanese map
$$j^{et}_{loc}:x \mapsto [P^{et}(x)] \in H^1_f(G_v,U^{et}).$$
When we let $L=F$, the number field itself
and let $b\in \cX(R)$ be an integral point,
then the $\G=\Gal (\bar{F}/F)$-action factors through
$G=G_T$ (as in \cite{wildeshaus},
proof of theorem 2.8), with notation following the introduction.
As explained in loc. cit., the point here is that there is an exact sequence
$$0\ra U^{et} \ra U^{et, ar} \ra G_T \ra 0$$
that is split by the base-point $b$,
 accounting for the $\G$ action on $U^{et}$
and also for the twisted actions corresponding to other points.
(Here, the middle term is an `arithmetic' fundamental group
obtained by pushing out a pro-finite arithmetic fundamental group $\hat{\pi}_1(\cX,b)$
using $U^{et}$.)
Therefore, we have a global continuous cohomology
pro-variety
$$H^1(G_T, U^{et})$$
and a classifying sub-variety
$$H^1_f(G_T,U^{et})\subset H^1(G_T, U^{et}) ,$$defined to be the inverse image of
$H^1_f(G_v, U^{et})$ under the restriction map
$$H^1(G_T,U^{et})\ra H^1(G_v, U^{et}).$$
 Using this environment, we can finally define the global unipotent
Albanese map
$$j^{et}_{glob}:x\in \cX(R) \mapsto [P^{et}(x)] \in H^1_f(G_T,U^{et}).$$

To conclude, we have constructed a commutative diagram

$$\begin{array}{ccccc}
\cX(R)& \ra &\cX_v(R_v) & \ra & Res^{F_v}_{\Q_p}(U^{DR}/F^0) (\Q_p)\\
\da&  &\da & \nearrow & \\
H^1(G_T,U^{et})& \ra & H^1_f(G_v, U^{et})(\Q_p) & &
\end{array}$$
where the commutativity of the triangle on the right is in
\cite{faltings} and
\cite{olsson}, following earlier cases of \cite{shiho},
\cite{vologodsky}.

Passing from here to the  finite-level quotients for each $n$
gives us the diagram from the introduction:
$$\begin{array}{ccccc}
\cX(R)& \ra &\cX_v(R_v) & \ra & Res^{F_v}_{\Q_p}(U^{DR}_n/F^0) (\Q_p)\\
\da&  &\da & \nearrow & \\
H^1(G_T,U_n^{et})& \ra & H^1(G_v, U_n^{et})(\Q_p) & &
\end{array}$$
\section{Comments I}
The proof of the statement that conjecture 1 implies the
finiteness of $\cX(R)$
is exactly as in \cite{kim} and in Chabauty's original argument.
The dimension assumption implies that for $n$ as in the conjecture,
the image of the map
$$H^1_f(G_T, U^{et}_n) \ra U^{DR}_n/F^0$$
is contained in a proper Zariski-closed subset.
Hence, so is the image of
$$\cX(R) \ra  U^{DR}_n/F^0.$$
We deduce that there is a non-zero algebraic function on
$U^{et}_n/F^0$ that vanishes on this image.
But by theorem 1, this function is not identically
zero on any  residue disk of
$\cX(R_v)$ and is a Coleman function \cite{besser}. Therefore,
its zero set on $\cX(R_v)$ is finite.
It is perhaps worth noting that
$\cX(R_v) \ra U^{DR}_n/F^0$ is at most finite-to-one
(this follows from the $n=2$ case).
The essential point we have shown then is
that
$$\mbox{Im} (\cX(R)) \subset
\mbox{Im} (\cX(R_v))\cap \mbox{Im}(H^1_f(G_T, U^{et}_n))$$
and that the latter intersection is finite (of course assuming conjecture 1).
Therefore, the strategy of this paper
is reminiscent of  Serge Lang's suggestion \cite{lang}
to prove directly that the intersection of a curve with the Mordell-Weil
group of its Jacobian is finite, except that the approach here is
p-adic and non-abelian.

We proceed to comment on the
relation between conjecture 1 and various
standard conjectures.
Such a connection is based upon a study of various
groups of Selmer type. Given a $\Q_p$-representation $V$ of
$G_{T}$ we will be especially interested in
$$Sel^0_T(V):=\mbox{Ker}(H^1(G_T, V) \ra \oplus_{w\in T} H^1(G_w, V))$$
\begin{conj}
Let $X/F$ be a smooth curve and
let $V_n=H^1_{et}(\bar{X}, \Q_p)^{\otimes n}(1)$. Then
$Sel^0_T(V_n)=0$ for $n>>0$.
\end{conj}
This conjecture appears to be a special case of some
general expectations about motivic representations,
and I am not sure to whom it should be attributed.
Notice that
$V_n$ is a direct summand of $H^n(\bX^n, \Q_p(1))$.
There is a Chern class map from $K$-theory
$$ch_{n,r}: K^{(r)}_{2r-n-1}(X^n)\ra H^1(\Gal(\bar{F}/F), H^n(\bX^n, \Q_p(r)))$$
for $2r-1\neq n$
and cycle maps
$$c_r: CH^r_0(X^n) \ra H^1(\Gal(\bar{F}/F), H^{2r-1}(\bX^n, \Q_p(r)))$$
where the superscript in the $K$-group
refers to an associated graded space  for the gamma filtration, while the
subscript in the Chow group refers to the cycle classes homologically
equivalent to zero.
The images of $ch_{n,r}$ and $c_r$ lie inside the subspaces
$H^1_g(\Gal(\bar{F}/F),H^n(\bX^n, \Q_p(r))) $ consisting of cohomology classes
that are unramified at almost all primes and potentially semi-stable at
all primes. In particular, $$Sel^0_T(V_n)\subset H^1_g(\Gal(\bar{F}/F),H^n(\bX^n, \Q_p(1))).$$

In the course of formulating their famous conjectures,
Bloch and Kato conjectured (\cite{B-K}, conjecture 5.3) that $ch_{n,r}$ induces an isomorphism
$$K^{(r)}_{2r-n-1}(X^n)\otimes \Q_p \simeq H^1_g(\Gal(\bar{F}/F),H^n(\bX^n, \Q_p(r)))$$
This can be view as a grand generalization of the
finiteness conjecture for Tate-Shafarevich groups of Abelian varieties.
It is also stating that certain extensions are geometric in a
specific sense.
We will explain below how the section conjecture can also be viewed
in a similar light. In any case,
if $n>2r-1$ the $K$-groups are zero.
Therefore
\begin{obs}
The conjecture of Bloch and Kato implies conjecture 2.
\end{obs}

Tony Scholl has pointed out to me that
\begin{obs}
Conjecture 2 is also implied by the Fontaine-Mazur conjecture.
\end{obs}
This is because an extension
$$0\ra V_n \ra E \ra \Q_p \ra 0$$
corresponding to an element of this group must be
semi-stable at $p$, and hence, a $\Q_p$ linear combination
of extensions of geometric origin. In particular, it must
carry a weight filtration. $V_n$ is pure of weight
$n-2$. Since a   map between mixed representations is strict for
the weight filtration, there exists a vector $v\in W_0E$
that maps to $1\in \Q_p$. If $n\geq 3$, then
$W_0E\cap V_n=0$, and hence, $W_0E$ must be a
1-dimensional complement to $V_n$.
Thus, the extension must split.

We will show
\begin{prop}
Conjecture 2 implies conjecture 1.
\end{prop}

If we use the symbol $U^M$ for
the fundamental group in any of the realizations of interest  and
$U^M_n=Z^{n}\backslash U^M$, then we have an exact sequence
$$0 \ra Z^{n+1}\backslash Z^n \ra U^M_{n+1} \ra U^M_{n}\ra 0$$
We need to calculate the dimension $d_n$ of  $Z^{n+1}\backslash Z^n$.
This is achieved as follows: $U^M$ is the unipotent completion of either
a free group on $m=2g+s-1$ generators, where $g$ is the genus of $X'$
and $s$ the order of $X'\setminus X$, or a free group on
$2g$ generators modulo a single relation of degree 2 (the compact case).
According to  \cite{labute} $d_n$ is given by the recursive formula
$$\Sigma_{k|n}kd_k=m^n$$
in the open case and
$$\Sigma_{k|n}kd_k=(g+\sqrt{g^2-1})^n+(g-\sqrt{g^2-1})^n $$
in the compact case. So one gets the asymptotics
$d_n \approx m^n/n$ in the non-compact case and
$d_n\approx (g+\sqrt{g^2-1})^n/n$ in the compact case.
(The small difference between the two formulas appears to be significant.
Even after considerable effort, the argument of observation 3 below could not
be adapted to the compact case.)
Notice here that in the compact case, $g\geq 2$.
We can also estimate $F^0(Z^{n+1}\backslash Z^n)$ in the De Rham realization
by noting that it is a quotient of
$H_1^{DR}(X_v)^{\otimes n}$ which has the filtration
dual to that on $H^1_{DR}(X_v)^{\otimes n}$. So $F^0 H_1^{DR}(X_v)^{\otimes n}$
has dimension equal to
 the codimension of
$$F^1H^1_{DR}(X_v)^{\otimes n}$$
$$=[F^1H^1_{DR}\otimes (H^1_{DR})^{\otimes (n-1)}]
\oplus [H^1_{DR}\otimes F^1H^1_{DR} \otimes (H^1_{DR})^{\otimes (n-2)}]
\oplus \cdots \oplus [(H^1_{DR})^{\otimes (n-1)}\otimes F^1H^1_{DR}]$$
which therefore is equal to
$$\dim H^1(\O_{X'_v})^{\otimes n}=g^n.$$
Therefore, $F^0(Z^{n+1}\backslash Z^n)\leq g^n$. In particular, it
does not contribute
to the asymptotics. That is, the jump in dimensions as one goes
from $U^{DR}_{n}$ to $U^{DR}_{n+1}$ is determined by the asymptotics of $d_n$.
On the other hand, in the \'etale realization, we have
$$0 \ra H^1(G_T,Z^{n+1}\backslash Z^n) \ra H^1(G_T, U^{et}_{n+1}) \ra H^1 (G_T, U^{et}_{n})$$
in the sense of \cite{kim} section 1, whereby the middle term
is an $H^1(G_T,Z^{n+1}\backslash Z^n)$-torsor over a subvariety
of the last term.
So we can control the change in dimensions using the
 Euler characteristic formula
$$\dim H^1(G_T,Z^{n+1}\backslash Z^n)-\dim H^2(G_T,Z^{n+1}\backslash Z^n)=\dim
(Z^{n+1}\backslash Z^n)^{-},$$
where the negative superscript refers to the (-1) eigenspace  of complex conjugation.
By comparison with complex Hodge theory, we see that the right hand side
is $d_n/2$ for $n$ odd. So it remains to observe that
\begin{lem} Conjecture 2 implies
$$\dim H^2(G_T,Z^{n+1}\backslash Z^n) \leq  P(n)g^n$$ for an effective
polynomial $P(n)$ of $n$.
\end{lem}

{\em Proof of lemma.}
We have a surjection
$$H^2(G_T, H_1(\bar{X}, \Q_p)^{\otimes n}) \ra H^2(G_T,Z^{n+1}\backslash Z^n)  \ra 0,$$
so it suffices to prove the estimate for the first group.
But by conjecture 2 and Poitou-Tate duality, we get an injection
$$H^2(G_T,  H_1(\bar{X}, \Q_p)^{\otimes n}) \hra \oplus_{w\in T} H^2(G_w,H_1(\bar{X}, \Q_p)^{\otimes n})$$
for $n$ sufficiently large.
Hence, we need only estimate the dimension of the local groups.
Applying local duality allows us to reduce this to
the study of
$$H^0(G_w, H^1(\bar{X}, \Q_p)^{\otimes n}(1)).$$
For $w$ not dividing $p$, a combination of the
geometric weight filtration and the monodromy filtration gives us

$$W: \ \ \ W_0 \subset W_1 \subset W_2=H^1(\bar{X}, \Q_p)$$
where
$\dim Gr^W_0\leq g$.
>From this, we get immediately that
$$\dim H^0(G_w, H^1(\bar{X}, \Q_p)^{\otimes n}(1)) \leq n \dim Gr^W_2
g^{n-1}+ \frac{n(n-1)}{2} (\dim Gr^W_1)^2 g^{n-2}
$$
For $w$ dividing $p$, we use the Hodge-Tate decomposition:
$$H^1(\bar{X},\Q_p)\otimes \C_p \simeq H^0(X', \Om_{X'}(\log D))\otimes \C_p (-1)\oplus H^1(X', \O_{X'})\otimes \C_p$$
as $G_w$-representations. Thus,
$$H^0(G_w, H^1(\bar{X}, \Q_p)^{\otimes n}(1)) \hra H^0(G_w,[H^0(X', \Om_{X'}(\log D))
\otimes \C_p (-1)\oplus H^1(X', \O_{X'})\otimes \C_p]^{\otimes n}(1))$$
from which we easily extract the estimate
$$\dim H^0(G_w, H^1(\bar{X}, \Q_p)^{\otimes n}(1)) \leq n (g+s-1)g^{n-1}$$
$\Box$

Since we have already discussed the structure of $U^M$
in the proof of the previous proposition,
it becomes easy to outline the relationship to
a conjecture of Jannsen. There, one starts with any smooth projective
variety $V$ over $F$ having good reduction outside
a finite set $T$ of primes (which we take to include all primes
dividing $p$)
and considers
$H^2(G_T, H^n(\bar{V}, \Q_p)(r))$
with various twists $r$.
Then Jannsen conjectures:

\begin{em}
$$H^2(G_T, H^n(\bar{V}, \Q_p)(r))=0$$
for $r\geq n+2$.
\end{em}

In fact, we need only a weaker variant

\begin{em} (Weak Jannsen conjecture) There exists a $k>0$
such that for varieties $V$ as above
$$H^2(G_T, H^n(\bar{V}, \Q_p)(r))=0$$
for $r\geq n+k$.
\end{em}

\begin{obs}
When $X$ is affine, Jannsen's conjecture in the weaker version
implies conjecture 1.
\end{obs}

The  point here is to try once more to estimate
the $H^2$. But in the affine case, the growth rate
of $Z^{n+1}\backslash Z^n$ is like $m^n/n$, so it will
certainly suffice to show that
$H^2(G_T, H_1(\bar{X}, \Q_p)^{\otimes n})$
is bounded by $Ca^n$ for some $a<m$.
Now, easy homological arguments give us the existence
of a constant $K$ such that
$$\dim H^2(G_T, V)\leq K \dim V$$
Recall that $m=2g+s-1$ and that
$H_1(\bar{X},\Q_p)$ fits into an exact sequence
$$0\ra W \ra H_1(\bar{X},\Q_p) \ra V \ra 0$$
where $W=\Q_p(1)^{s-1}$  and
$V=H_1(\bar{X'}, \Q_p)$. (Here, we have passed to a finite extension
where all the points at infinity are rational. This passage is
admissible for our $H^2$ considerations.)
Therefore, in the Grothendieck group,
we get
$$H_1(\bar{X},\Q_p)^{\otimes n}=V^{\otimes n} \oplus
n(V^{\otimes (n-1)}\otimes W )\oplus \cdots \oplus
{n\choose k}(V^{\otimes (n-k)}\otimes
W^k) \oplus \cdots$$
But $V\simeq H^1(\bar{X}', \Q_p)(1)$
so that
$V^{\otimes l}$ is a direct summand of $H^l((\bX')^l,\Q_p)(l)$.
Therefore, $H^2(G_T, V^{\otimes l}(r))=0$ for
$r\geq k$.
Hence,
$$\dim H^2(G_T, H_1(\bar{X},\Q_p)^{\otimes n})\leq P(n)
\Sigma_{l=0}^{k-1}\dim H^2(G_T, V^{\otimes (n-l)}(l))$$
for some polynomial $P(n)$.
We analyze the cohomology groups on the right.
 The Weil pairing induces an inclusion
$$\Q_p(1) \hra V^{\otimes 2}$$
which continues to an inclusion
$$\Q_p(k) \hra V^{\otimes (2k)}$$
with image  a direct summand.
Write
$V^{\otimes 2k}=M\oplus \Q_p(k)$.
One the one hand,  since any $M^{\otimes i}$ is a direct summand of a $V^{\otimes l}$, we have
$$H^2(G_T, V^{\otimes i}\otimes M^{\otimes j}(r))=0$$
for $r\geq k$. On the other,
$$V^{\otimes 2kn}=M^{\otimes n}\oplus M'(k)$$
where $M'$ is a direct sum of (non-negative) tensor products of $M$ and $\Q_p(1)$.
Hence,
$$V^{\otimes (2kn+i)}=M^{\otimes n}\otimes V^{\otimes i}\oplus
M'\otimes V^{\otimes i} (k)$$
and
$$H^2(G_T, V^{\otimes (2kn+i)}(l))=H^2(G_T, M^{\otimes n}\otimes V^{\otimes i}(l))$$
for $l\geq 0$. Therefore, as $i$ runs through
$0, \ldots, 2k-1$, there is a  polynomial $Q$ such that
$$\dim H^2(G_T, V^{\otimes (2kn+i)}(l))\leq Q(n)
 \dim M^{\otimes n}=Q(n)((2g)^{2k}-1)^n\leq Q(n)((2g)^{2k}-1)^{(2kn+i)/(2k)}.$$
Rewriting the exponent, for any $n\geq 0$,
$$\dim H^2(G_T, V^{\otimes n}(l))\leq Q(n) [((2g)^{2k}-1)^{1/(2k)}]^n$$
from which the desired estimate follows.

\section{Comments II}
It is rather interesting to consider the special case
of an elliptic curve minus the origin.
For example, our method shows  finiteness
of integral points for CM-elliptic curves $E$ defined
over $\Q$ of rank 1, independently of any conjecture.
This follows from analyzing somewhat carefully
a refined Selmer variety
$$Sel_T(U_n),$$ defined to be the classes in
$H^1_f(G_T, U^{et}_n)$ that are potentially unramified
at primes not equal to $p$ and potentially
crystalline at $p$. (Here of course, the $U$ refer to the unipotent
fundamental groups of $E$ minus the origin.)
The definition of this subvariety depends
on  the construction of certain
local subvarieties
$$H^1_{pf}(G_v, H_v;  U^{et}_n)\subset H^1(G_v, U^{et}_n)$$
($v|p$) and
$$H^1_{pun}(G_v, H_v;U^{et}_n)\subset H^1(G_v, U^{et}_n)$$
($v$ not dividing $ p$) for various
$v$ and subgroups of finite index
 $H_v\subset G_v$. Since the pattern is exactly as in \cite{kim}, section 1,
we will omit the details. The point is
that if $N$ is a closed subgroup of $G$, then
the restriction map
$H^1(G, U^{et}_n) \ra H^1(N,U^{et}_n)$
is functorial.
We can therefore construct the unramified cohomology
varieties $H^1_{un}(G_v, U^{et}_n)$ as the inverse image
of the base-point in $H^1(I_v, U^{et}_n)$, where $I_v\subset G_v$
is the inertia subgroup. For this last object, we need only
deal with it when the action of $I_v$ on $U^{et}_n$ is trivial.
In this case, the cohomology functor  is not necessarily representable
by a variety. However, we still have functorial  sequences
$$0\ra H^1(I_v, Z^{n+1}\bs Z^n) \ra H^1(I_v, U^{et}_{n+1}) \ra H^1(I_v, U^{et}_n)$$
which are exact in the naive sense that the inverse image of the base point
under the second map is the kernel. (However, the other fibers may be empty.)
Now, the condition that a class $c\in H^1(G, U^{et}_{n+1})$ goes to
zero in $H^1(I_v, U^{et}_{n+1})$ is an intersection of the condition
that it goes to zero in $H^1(I_v, U^{et}_n)$, which defines a subvariety
by induction, and then the condition that the image in
$H^1(I_v, Z^{n+1}\bs Z^n)$ is trivial. This last space is a vector group,
and hence, the resulting condition defines the desired `unramified'
subvariety at level $n+1$.

Then, by considering the restriction maps
$$r_v:H^1(G_v, U^{et}_n) \ra H^1(H_v, U^{et}_n),$$
 to subgroups of finite index, we get the potentially unramified subvarieties
$$H^1_{pun}(G_v, H_v; U^{et}_n):=
r_v^{-1}(H^1_{un}(H_v, U^{et}_n))$$ and  potentially crystalline
subvarieties $$H^1_{pf}(G_v, H_v; U^{et}_n):=
r_v^{-1}(H^1_{f}(H_v, U^{et}_n)).$$
We know that
$E$ acquires good reduction everywhere after base extension
to the field $F$ generated by the 12-torsion. We then take
the $H_v$ to be the Galois groups of the completions
$F_v$ for each place $v$ of $F$ lying above a place of $T$.
Finally, we define
$$H^1_{pf}(G_p, U^{et}_n):=\cap_{v|p} H^1_{pf}(G_p, H_v, U^{et}_n),$$
and
$$H^1_{pun}(G_v, U^{et}_n):=\cap_{w|v} H^1_{pun}(G_v, H_w, U^{et}_n),$$
for $v \neq p$. Then the global group
$Sel_T(U_n)$
is defined to be the intersection of the inverse images
of all these local groups under the global-to-local
 restriction maps:
$$Sel_T(U_n):=\cap_{v\in T}[(Res^{glob}_{loc,v})^{-1}(H^1_{pun/ pf}(G_v, U^{et}_n))]\subset H^1(G_T,U^{et}_n).$$
With these definitions in place we claim
there is an exact sequence
$$0\ra Sel_T(Z^{n+1}\backslash Z^n) \ra Sel_T(G_T, U^{et}_{n+1})
\ra Sel_T(U^{et}_{n})$$
again in the sense that the vector group on the
left acts freely on the middle with quotient
a subvariety of the right hand side.
Most of the proof is a consequence of the exactness of
$$0\ra H^1(G_T, Z^{n+1}\backslash Z^n) \ra H^1(G_T, U^{et}_{n+1}) \ra H^1(G_T, U^{et}_{n}),$$
(cite{kim}, section 1)
leaving us to check that if
$c_1, c_2$ are two elements
of $Sel_T(G_T, U^{et}_{n+1})$ mapping to the same class
in $Sel_T(G_T, U^{et}_{n})$, then one is a translate of
the other by a class in $Sel_T(Z^{n+1}\backslash Z^n)$.
This can be checked locally at each prime.
So let $w$ be a prime of $F$ lying above
$p$ and let $c_1,c_2 \in Z^1(H_w, U^{et}_{n+1})$ be
two cocycles representing classes of $H^1_f(H_w, U^{et}_{n+1})$
and assume they have the same image in $H^1_f(H_w, U^{et}_{n})$.
We already know that $c_2=c_1z$ for some cocycle
$z\in Z^1(H_w, Z^{n+1}\backslash Z^n)$ and we need to check
that it trivializes over $B_{cr}$.
Since both $c_1$ and $c_2$ trivialize over
$B_{cr}$, there are elements $u_1, u_2 \in U^{et}_{n+1}(B_{cr})$
such that $c_i(g)=g(u_i^{-1})u_i$ for $g\in H_w$.
Then $z(g)=g(u_1u_2^{-1})u_2u_1^{-1}$.
Therefore, we see that $u_2u_1^{-1}$ maps to an $H_w$-invariant
element in $U^{et}_{n}(B_{cr})$.
But the invariant points just comprise $U^{DR}_{n}$ and the map
$U^{DR}_{n+1}\ra U^{DR}_{n}$
is surjective. Therefore, there is
an $H_w$-invariant element $v\in U_{n+1}(B_{cr})$
such that $u_2u_1^{-1}=vc$ for some $c\in Z^{n+1}\backslash Z^n (B_{cr})$.
Hence,  $z(g)=g((vc)^{-1})(vc)=g(c^{-1})c$, giving us the triviality
we want. (Of course we should have included a $\Q_p$-algebra $R$
in the notation to prove all this functorially, but we have omitted
it for brevity.)
At a prime $w$ not dividing $p$, the condition to check concerns
unramified cocycles, and hence is easier:
Suppose $c_1, c_2$ are two  unramified
cocycles with values in $U^{et}_{n+1}$ that map to
the same class in $H^1(H_w, U^{et}_{n})$.
Then $c_2=c_1z$ for some cocycle $z$ with values in $Z^{n+1}\backslash Z^n$.
But since the action on $U^{et}_{n+1}$ itself is
unramified, all cocycles restricted to
 the inertia group $I_w$ are homomorphisms.  Then since $c_1$ and $c_2$
are in fact trivial on $I_w$, so is
$z$.

For $n=2$, we get from this the inequality of dimensions
$$\dim Sel_T( U^{et}_3)\leq \dim Sel_T( U^{et}_2)
+\dim Sel_T(Z^3\backslash Z^2).$$
However, $$Z^3\backslash Z^2\simeq H_1(\bar{E},\Q_p)\wedge H_1(\bar{E},\Q_p)
\simeq \Q_p(1).$$
This implies that
$ Sel_T(Z^{3}\backslash Z^2)$ consists of the units in $\Z$
tensor $\Q_p$, and hence, is zero (\cite{B-K}, example 3.9).
Therefore,
$$\dim Sel_T(U_3^{et})=\dim Sel_T(U_2^{et})=
\dim Sel_T(H_1(\bar{E}, \Q_p))=1,$$ by the rank one hypothesis.

On the De Rham side, we have the same dimension count $(=1)$
for $Z^3\backslash Z^2$. Meanwhile, considering the
quotient
map $$H_1^{DR}(E_p)\otimes H_1^{DR}(E_p) \ra Z^3\backslash Z^2$$
and the fact that $F^0H_1\wedge F^0H_1=0$,
we get $F^0(Z^3\backslash Z^2)=0$. So
$$\dim F^0U^{DR}_3=\dim F^0U^{DR}_2=1$$ and
$$\dim U^{DR}_3/F^0=3-1=2.$$
Therefore, $j^{DR}_3$ suffices to finish the job.

This calculation shows that  in explicit situations like this with a
fixed $n$, it is worth imposing
extra Selmer conditions to improve the dimension estimate, in contrast to
the general finiteness problem where refined conditions
 appear to be irrelevant.
It is probably also worth looking into this case more carefully
with a view towards working out explicit bounds in the
manner of Coleman's refinement of Chabauty \cite{coleman}. The suspicion is that this
will require understanding the Dieudonn\'e map $D$ more precisely,
along the lines of `non-abelian
explicit reciprocity laws' \`a l\`a Kato \cite{kato}.

It was pointed out to me by John Coates that
the general CM case is intimately related
to a pseudo-nullity conjecture from the Iwasawa
theory of elliptic curves. That is, let
$F_{\infty}$ be the field generated by the $p$-power torsion
of $E$, $F_1$ the field generated by the $p$-torsion,
$G_{\infty}=\Gal(F_{\infty}/F_1)$,
and $\Lambda=\Z_p[[ G_{\infty}]]$ the Iwasawa algebra
of $G_{\infty}$.
Let $I$ be  the Galois group of the maximal abelian unramified pro-p extension
of $F_{\infty}$ split over the primes lying above $T$.
Coates and Sujatha have conjectured \cite{C-S}
that $I$ is pseudo-null as a $\Lambda$-module.

\begin{obs}
This pseudo-nullity implies conjecture 1 for $E\setminus \{0\}$.
\end{obs}

Recall that what is necessary is to estimate
$H^2(G_T, Z^{n+1}\backslash Z^n)$
which is implied by an estimate for
$H^2(G_T, H_1(\bar{E}, \Q_p)^{\otimes n})$. But for just this part,
we can assume that the base-field is $F_1$
since the corestriction map is surjective on $H^2$.
In fact, the pseudo-nullity implies
that
$$\dim H^2(G_T, H_1(\bar{E}, \Q_p)^{\otimes n})\leq P(n)$$
for some polynomial $P(n)$.
To see this, first note that
it suffices to prove such an inequality for the Selmer group
$$Ker( H^2(G_T, H_1(\bar{E}, \Q_p)^{\otimes n})\ra \oplus_{v\in T}
H^2(G_v, H_1(\bar{E}, \Q_p)^{\otimes n}))$$
since the local contributions can be bounded exactly
as in section 3 (note that $g=1$ in this case).
So it suffices to bound the $H^1$ Selmer group
$$Sel^0_T(H^1(\bar{E}, \Q_p)^{\otimes n}(1))=Ker( H^1(G_T, H^1(\bar{E}, \Q_p)^{\otimes n}(1))\ra \oplus_{v\in T}
H^1(G_v, H^1(\bar{E}, \Q_p)^{\otimes n}(1))),$$
and we have an isomorphism:
$$Sel^0_T(H^1(\bar{E}, \Q_p)^{\otimes n}(1))\simeq \Hom_{\Lambda}
(I, H^1(\bar{E}, \Q_p)^{\otimes n})(1))$$
coming from the Hochschild-Serre spectral sequence.
To see this, it is perhaps convenient to
utilize the language of continuous \'etale cohomology \cite{jannsen2}
applied to the inverse system
$$(T_i)=(H^1(\bar{E}, \Z/p^i)^{\otimes n}(1)).$$
The Hochschild-Serre spectral sequence in this context gives us
the exact sequence
$$0\ra H^1(G_{\infty}, (H^0(N,T_i))) \ra H^1(G_T, (T_i)) \ra H^0(G_{\infty}, (H^1(N, T_i))) \ra
H^2(G_{\infty}, (H^0(N,T_i)))$$
(\cite{jannsen2}, theorem 3.3), where $N\subset G_T$
is the subgroup fixing $F_{\infty}$. Since the $N$-action is trivial,
all the $H^0(N,T_i)$ are just $T_i$, and hence, satisfy
the Mittag-Leffler condition. Therefore, all the terms but one  become readily
identified with the continuous Galois cohomology of $T=\invlim T_i=H^1(\bar{E}, \Z_p)^{\otimes n}(1)$
while
$$H^0(G_{\infty}, (H^1(N, T_i))):=\invlim \Hom_{G_{\infty}}(N, T_i)=\Hom_{G_{\infty}}(N, T).$$
When we tensor with $\Q$, we get therefore
$$H^1(G_T,H^1(\bar{E}, \Q_p)^{\otimes n}(1))\simeq \Hom_{\Lambda}(N, H^1(\bar{E}, \Q_p)^{\otimes n}(1)).$$
Passing  to the Selmer group,
the last group gets replaced by
$\Hom_{\Lambda}
(I, H^1(\bar{E}, \Q_p)^{\otimes n}(1))$.

But the pseudo-nullity implies that
$I$ is finitely generated over $\Z_p$ (see, e.g., \cite{greenberg}, section 2, paragraph following
lemma 2).
Thus, if we fix
a Frobenius element $F$ in
$G_{\infty}$, then only finitely many eigenvalues
can occur for $F$ in $I\otimes \Q_p$.
Weight considerations imply then that
$$ \Hom_{\Lambda}
(I, H^1(\bar{E}, \Q_p)^{\otimes n}(1))=0$$
for $n>>0$.

\section{Comments III}
Let us consider the section conjecture briefly in terms of Galois cohomology in order
to relate it to the Selmer varieties that occur in this paper.
For simplicity, we consider it in the situation where
a point $b$ exists in $X(F)$ and is used as the base point for
the fundamental group.
So we have the exact sequence
$$0\ra \hat{\pi}_1(\bX,b) \ra \hat{ \pi}_1(X,b) \ra \G \ra 0$$
of profinite fundamental groups from the introduction
as well as one fixed splitting $$s:\G \ra \hpi_1(X,b).$$
We use this splitting to make
$\G$ act on $\hat{\pi}_1(\bX,b) $ by conjugation via $s$.
Suppose $t:\G \ra \hat{\pi_1}(X,b)$ is another splitting.
Then for each $g \in \G$, there is a
$z_t(g) \in \pi_1(\bX,b)$ such that
$t(g)=z_t(g)s(g)$. We have
$t(g)t(h)=t(gh)=z_t(gh)s(gh)$
and hence, $z_t(gh)s(g)s(h)=z(g)s(g)z_t(h)s(h)$,
or $z_t(gh)=z_t(g)g(z_t(h))$.
That is to say, $z_t$ is a 1-cocycle.
In fact, a straightforward computation shows
that this map sets up a bijection
$$[Split(X)/\sim] \simeq H^1(\G, \hat{\pi}_1(\bX,b))$$
to a continuous cohomology set.
In the case of a point $x\in X(F)$,
the corresponding cohomology class is the class of
the torsor $\hat{\pi}_1(\bX;x,b)$ of \'etale paths
from $b$ to $x$. So we have simply reformulated
the section conjecture in terms
of
a classifying map
$$X(F) \ra H^1(\G, \hat{\pi}_1(\bX,b))$$
that goes to a classifying space for torsors
similar to the ones occurring in this paper.
In fact, this map is easily seen to be injective
so that the content of the section conjecture
lies in surjectivity. That is, it is saying that
all torsors are of `geometric origin'
in a very precise sense.

Let us try to formulate the correct analogue of the
section conjecture for an elliptic curve
$E$
equipped with the origin $o\in E$ (which we are {\em not} deleting).
We still have an exact sequence
$$0\ra \hpi_1(\bar{E}, o) \ra \hpi_1(E,o) \ra \G \ra 0$$
and the splittings are still classified by
$H^1(G, \hpi_1(\bar{E}, o))$.
However, in this non-hyperbolic situation, the splittings that come from
points satisfy local conditions
that amount to exactly the subspace
$$H^1_g(G, \hpi_1(\bar{E}, o))\subset H^1(G, \hpi_1(\bar{E}, o))$$
studied by Bloch-Kato. That is, the correct classifying map to
study in this case
is
$$E(F) \ra H ^1_g(G, \hpi_1(\bar{E}, o)).$$
 Since
$$\hpi_1(\bar{E}, o)=H_1(\bar{E},\hat{\Z})\simeq H^1(\bar{E},\hat{\Z})(1),$$
we conclude that the elliptic curve
analogue of the section conjecture
is nothing but the finiteness conjecture for the
Shafarevich group of $E$. Of course the map that occurs here
is identified with the cycle map coming from
$CH^1_0(E)\simeq E(F)$. In this abelian situation,
the Chern class maps from higher $K$-theory appear not to play a
role in the `section conjecture.'

In view of  the analogy just outlined,
it is perhaps reasonable to think of
the relevant   parts of the
Bloch-Kato conjecture in the hyperbolic case as
comprising a `linear analogue' of the section conjecture.
In other words, one possible  perspective
is that we have proved an implication
of the form

\medskip
`linear' section conjecture $\Rightarrow$ Faltings' theorem.
\medskip

Within this restricted context, one can consider the
Fontaine-Mazur conjecture as being a `weak' Bloch-Kato/BSD
conjecture, in that one is merely asking that certain
extensions be geometric in whatever manner, as opposed to
the precise manner of coming from points or the $K$-theory of
 specific varieties under discussion. Therefore,
 it seems natural
to ask for
a weak section conjecture, whereby
one demands that every class in
$H^1(\G, \hpi_1(\bX,b))$ be of geometric origin
in some definite sense that lacks as yet a good formulation.
(For example, coming to grips with such a formulation
 appears to be considerably more difficult
 then a well-known notion from the theory of mixed motives:
` a system of realizations
of geometric origin.')
Meanwhile, for these vague ideas to acquire more weight,
it would be nice to prove an implication of the form
\medskip

section conjecture $\Rightarrow$ linear section conjecture.
\medskip

\begin{flushleft}
{\bf Acknowledgements:}

-The author was supported in part by
a grant from the National Science Foundation.
\smallskip

Deep gratitude is due to:
\smallskip

-Shinichi Mochizuki for overarching
influence stemming from his `non-linear' philosophy;

-Dick Hain for continuing  to provide patient answers to questions
on Hodge theory;

-Bill McCallum, for emphasizing over the years the importance
of Chabauty's method;

-Dinesh Thakur, for a careful and enthusiastic reading of the manuscript;

-Martin Olsson, for communicating to me his deep results on
non-abelian $p$-adic Hodge theory;

-Uwe Jannsen, Guido Kings, Jan Nekovar, and Tony Scholl, for helpful discussions
on mixed motives;

-Gerd Faltings, Florian Pop,  John Coates, and Kazuya Kato
for persistent encouragement of this research and many enlightening
discussions on fundamental groups and Iwasawa theory;

-and, of course, Serge Lang, for too many things to be
listed.

\end{flushleft}

{\footnotesize DEPARTMENT OF MATHEMATICS, PURDUE UNIVERSITY,
  WEST LAFAYETTE, INDIANA 47907  and DEPARTMENT OF MATHEMATICS, UNIVERSITY OF ARIZONA,
TUCSON, AZ 85721, U.S.A. }

{\footnotesize EMAIL: kimm@math.purdue.edu}

\end{document}